# The rational abelianization of the Chillingworth subgroup of the mapping class group of a surface

RYOTARO KOSUGE

The Chillingworth subgroup of the mapping class group of a compact oriented surface of genus $g$ with one boundary component is defined as the subgroup whose elements preserve nonsingular vector fields on the surface up to homotopy. In this work, we determine the rational abelianization of the Chillingworth subgroup as a full mapping class group module. The abelianization is given by the first Johnson homomorphism and the Casson–Morita homomorphism for the Chillingworth subgroup. Additionally, we compute the order of the Euler class of a certain central extension related to the Chillingworth subgroup and determine the kernel of the Casson–Morita homomorphism for the Chillingworth subgroup.

## 1 Introduction

Throughout this paper, we assume that all surfaces are compact, connected, and oriented. Let $\Sigma_{g,1}$ (resp. $\Sigma_{g,*}$, $\Sigma_g$) denote a surface of genus $g$ with one boundary component (resp. with a fixed base point, or with no boundary and no fixed point). The *mapping class group*, denoted by $\mathcal{M}_{g,1}$ (resp. $\mathcal{M}_{g,*}$, $\mathcal{M}_g$), is defined as the group of isotopy classes of orientation-preserving self-diffeomorphisms of the surface that fix the boundary or the base point pointwise. For oriented surface bundles, the structure group is the orientation-preserving diffeomorphism group of the surface. Except for a finite number of cases where the genus is small, this diffeomorphism group is homotopy equivalent to the mapping class group, which is discrete. As a result, their classifying spaces are homotopy equivalent, and therefore, the group cohomology of the mapping class group is equivalent to the characteristic classes of surface bundles. The mapping class group naturally acts on various structures of the surface. For example, it acts on the first integral homology group of the surface $H=H_1(\Sigma_{g,1};\mathbb{Z})$, preserving the intersection form. Consequently, the mapping class group acts on $H$ via the integral symplectic group $\mathrm{Sp}(2g,\mathbb{Z})$, called the symplectic representation. Through this, several important modules of the mapping class group can be

described by using representations of the integral symplectic group $\mathrm{Sp}(2g,\mathbb{Z})$, or the rational symplectic group $\mathrm{Sp}(2g,\mathbb{Q})$ over $\mathbb{Q}$.

## 1A Chillingworth subgroups and related background

Chillingworth [7; 8] studied the action of the mapping class group on the set of homotopy classes of nonsingular vector fields, focusing on winding numbers. This action is described by a crossed homomorphism, known as the *Chillingworth homomorphism*, which maps to the first integral cohomology group of the surface. The *Chillingworth subgroup* is defined as the subgroup of the mapping class group consisting of elements that preserve vector fields on the surface up to homotopy (see Proposition 3.3 for alternative definitions). Chillingworth subgroups of $\mathcal{M}_{g,1}$, $\mathcal{M}_{g,*}$, and $\mathcal{M}_g$ are denoted by $\mathrm{Ch}_{g,1}$, $\mathrm{Ch}_{g,*}$, and $\mathrm{Ch}_g$, respectively. Here, $\mathrm{Ch}_{g,*}$ and $\mathrm{Ch}_g$ are defined via certain natural homomorphisms $\mathcal{M}_{g,1}\to\mathcal{M}_{g,*}$ and $\mathcal{M}_{g,*}\to\mathcal{M}_g$ between mapping class groups, where $\mathcal{M}_{g,1}\to\mathcal{M}_{g,*}$ is obtained by collapsing the boundary to a point, and $\mathcal{M}_{g,*}\to\mathcal{M}_g$ is obtained by forgetting the base point.

Johnson [17] discussed the kernel of the Chillingworth class on the Torelli group, where the Chillingworth class is defined as the Poincaré dual of the Chillingworth homomorphism. Trapp [41] introduced a $(2g+1)$-dimensional linear representation of the mapping class group $\mathcal{M}_{g,1}$, referred to as Trapp's representation. In that work, he used this representation to study the action of the mapping class group on the first homology group of the unit tangent bundle of the surface and characterized the Chillingworth subgroup as the kernel of this linear representation. Furthermore, the Chillingworth subgroup has been studied in other contexts. Childers [6] studied its relationship with the subgroup generated by the simply intersecting pair (SIP) maps. Blanchet, Palmer and Shaukat [5] mentioned it in the context of the action of the mapping class group on the Heisenberg group of the surface, which is defined as a certain quotient of the surface braid group or a certain central extension of the first integral homology group of the surface by the infinite cyclic group. However, the structure of the Chillingworth subgroup has not been well studied.

Before we get into the main topic of this paper, we will introduce some background information. The Chillingworth subgroup is an intermediate-sized group between two significant subgroups in the context of the mapping class group: the *Torelli group* $\mathcal{I}_{g,1}$ which is defined as the kernel of the action of the mapping class group on the first homology group of the surface and the *Johnson kernel* $\mathcal{K}_{g,1}$ which is defined as the subgroup generated by Dehn twists along separating simple closed curves on the surface. Specifically, $\mathcal{K}_{g,1}\subset\mathrm{Ch}_{g,1}\subset\mathcal{I}_{g,1}$.

The structure of the rational abelianization of the Torelli group as a mapping class group module was determined by Johnson [19] using the *Johnson homomorphism* $\tau_{g,1}(1):\mathcal{I}_{g,1}\to\bigwedge^3 H_1(\Sigma_{g,1};\mathbb{Z})$, which he introduced and is now known as the first Johnson homomorphism. The target space is the third exterior power of the first homology group of the surface with the mapping class group acting naturally on it as the symplectic group. Since the Torelli group is a normal subgroup of the mapping class group, the mapping class group acts on it by conjugation. Under these actions, the first Johnson homomorphism is equivariant with respect to the action of the mapping class group. The structure of the rational abelianization of the Johnson kernel as a mapping class group module was determined by Dimca–Hain–Papadima [9], Morita–Sakasai–Suzuki [34] (in the case of closed surfaces without a base point), and Faes–Massuyeau [12, Theorem 3.2] (in the case of surfaces with one boundary component). The structure of the mapping class group module in this case is more complex but can be described as an extension of representations of the symplectic group (see Section 6).

The proof of Theorem A, which determines the rational abelianization of the Chillingworth subgroup, primarily involves analyzing the long exact sequence (inflation-restriction exact sequence)

$$H_2(\mathrm{Ch}_{g,1};\mathbb{Q}) \to H_2(U;\mathbb{Q}) \to H_1(\mathcal{K}_{g,1};\mathbb{Q})_U \to H_1(\mathrm{Ch}_{g,1};\mathbb{Q}) \to H_1(U;\mathbb{Q}) \to 0$$

induced by the first Johnson homomorphism $\tau_{g,1}(1):\mathrm{Ch}_{g,1} \to U$ restricted to the Chillingworth subgroup, where $U$ is the image $\tau_{g,1}(1)(\mathrm{Ch}_{g,1}) \subset \bigwedge^3 H_1(\Sigma_{g,1};\mathbb{Z})$.

We determine the rational abelianization of the Chillingworth subgroup using a result analogous to that of Hain [15] and the rational abelianization of the Johnson kernel by Faes and Massuyeau [12]. The former (Theorem B) corresponds to analyzing the leftmost map in the long exact sequence, while the latter corresponds to analyzing the third module from the left. Specifically, these are described by the (first) Johnson homomorphism $\tau_{g,1}(1)$ and the *Casson–Morita homomorphism* $d$. The first Johnson homomorphism is particularly important in the context of the Torelli group, as mentioned above, (see Section 2A); for example, the Johnson homomorphism for the Torelli group induces the rational abelianization of the Torelli group. The Casson–Morita homomorphism is closely related to the Casson invariant for homology 3-spheres (see Section 5) and provides one of the $\mathcal{M}_{g,1}$-invariant parts of the rational abelianization of the Johnson kernel $\mathcal{K}_{g,1}$.

In relation to the Casson–Morita homomorphism $d$, its properties on the Chillingworth subgroup, which are used in Theorems A and D, include its invariance under the action of the mapping class group and the determination of its image. These fundamental properties are summarized in Theorem C. Furthermore, although not directly relevant to Theorems A and D, Theorem C also includes an explicit description of the kernel of $d$ for the Chillingworth subgroup, as part of the fundamental properties of $d$.

Before presenting the theorems, we introduce some notation: $[-]_{\mathrm{Sp}}$ represents the linear representations of the rational symplectic group $\mathrm{Sp}(2g,\mathbb{Q})$ corresponding to Young diagrams. For details, see Section 4.

The rational abelianization of the Chillingworth subgroup of the mapping class group is as follows.

**Theorem A** *For $g \geq 6$, the rational abelianizations of the Chillingworth subgroups of the mapping class groups of the surfaces are induced by the Johnson homomorphisms and the Casson–Morita homomorphism*

$$\begin{aligned} d \oplus \tau_{g,1}(1) : \mathrm{Ch}_{g,1} &\to (\mathbb{Z} \oplus U) \otimes \mathbb{Q} \cong [0]_{\mathrm{Sp}} \oplus [1^3]_{\mathrm{Sp}}, \\ \tau_{g,*}(1) : \mathrm{Ch}_{g,*} &\to U \otimes \mathbb{Q} \cong [1^3]_{\mathrm{Sp}}, \\ \tau_g(1) : \mathrm{Ch}_g &\to \overline{U} \otimes \mathbb{Q} \cong [1^3]_{\mathrm{Sp}}, \end{aligned}$$

*where $U$ and $\overline{U}$ are images of the Chillingworth subgroups under the first Johnson homomorphisms. Specifically, their targets and the first rational homology groups of the Chillingworth subgroups are isomorphic as mapping class group modules.*

*In particular, the actions of the mapping class group on these abelianizations of the Chillingworth subgroups factor through the rational symplectic group* $\mathrm{Sp}(2g, \mathbb{Q})$*, and they decompose into irreducible representations of the rational symplectic group.*

**Theorem B** *The image (resp. kernel) of the homomorphisms between the second rational homology (resp. cohomology) induced by the first Johnson homomorphism*

$$\tau_{g,1}(1) = \tau_{g,1}(1)|_{\mathrm{Ch}_{g,1}} : \mathrm{Ch}_{g,1} \to U \subset \bigwedge^3 H_1(\Sigma_{g,1}; \mathbb{Z})$$

*for the Chillingworth subgroup for the genus-*$g$ *surface with one boundary is decomposed as mapping class group modules as*

$$\mathrm{Im}\big((\tau_{g,1}(1))_* : H_2(\mathrm{Ch}_{g,1}; \mathbb{Q}) \to H_2(U; \mathbb{Q})\big) = \begin{cases} [2^2 1^2]_{\mathrm{Sp}} \oplus [1^4]_{\mathrm{Sp}} \oplus [1^6]_{\mathrm{Sp}} & (g \geq 6), \\ [2^2 1^2]_{\mathrm{Sp}} \oplus [1^4]_{\mathrm{Sp}} & (g = 5), \\ [2^2 1^2]_{\mathrm{Sp}} & (g = 4), \\ \{0\} & (g = 3), \end{cases}$$

*and*

$$\mathrm{Ker}\big((\tau_{g,1}(1))^* : H^2(U; \mathbb{Q}) \to H^2(\mathrm{Ch}_{g,1}; \mathbb{Q})\big) = \begin{cases} [0]_{\mathrm{Sp}} \oplus [2^2]_{\mathrm{Sp}} \oplus [1^2]_{\mathrm{Sp}} & (g \geq 4), \\ [0]_{\mathrm{Sp}} \oplus [2^2]_{\mathrm{Sp}} & (g = 3). \end{cases}$$

*The same holds for the Chillingworth subgroup in the case of a fixed base point* $\mathrm{Ch}_{g,*}$*.*

Theorem B is used to prove Theorem A.

We examine the fundamental properties of the Casson–Morita homomorphism $d$ for the Chillingworth subgroup, focusing on explicitly determining its kernel, which is crucial to Theorem A.

**Theorem C** *The Casson–Morita homomorphism* $d = d|_{\mathrm{Ch}_{g,1}} : \mathrm{Ch}_{g,1} \to \mathbb{Z}$ *satisfies these properties:*

(1) *The Casson–Morita homomorphism* $d$ *is an* $\mathcal{M}_{g,1}$*-invariant homomorphism on the Chillingworth subgroup.*

(2) *The image* $\mathrm{Im}(d : \mathrm{Ch}_{g,1} \to \mathbb{Z})$ *of the Casson–Morita homomorphism for the Chillingworth subgroup corresponds to* $8\mathbb{Z}$*.*

(3) *The kernel* $\mathrm{Ker}(d : \mathrm{Ch}_{g,1} \to \mathbb{Z})$ *of the Casson–Morita homomorphism for the Chillingworth subgroup is given by the subgroup* $\langle T_{\gamma_1'} \rangle$ *generated by Dehn twists along the boundary of a genus-one subsurface with one boundary of the surface as shown in Figure 1, left, the normal subgroup* $\langle\langle B_0 \rangle\rangle \lhd \mathcal{M}_{g,1}$ *(recall that* $\langle\langle \bullet \rangle\rangle$ *denotes the normal closure) generated by a certain element* $B_0 := T_{\gamma_2'} T_{\gamma_3'}^{-1}$ *called the homological genus-zero bounding pair map as shown in Figure 1, right, and the commutator subgroup* $[\mathcal{K}_{g,1}, \mathcal{M}_{g,1}]$ *of the Johnson kernel and the full mapping class group as follows:*

$$\mathrm{Ker}(d : \mathrm{Ch}_{g,1} \to \mathbb{Z}) = \langle\langle B_0 \rangle\rangle \langle T_{\gamma_1'} \rangle [\mathcal{K}_{g,1}, \mathcal{M}_{g,1}].$$

Additionally, we compute the order of the Euler class of the natural central extension

$$0 \to \mathbb{Z} \to \mathrm{Ch}_{g,1} \to \mathrm{Ch}_{g,*} \to 1$$

related to the natural homomorphism $\mathcal{M}_{g,1} \to \mathcal{M}_{g,*}$. This is obtained by examining $d$ on the Chillingworth subgroup.

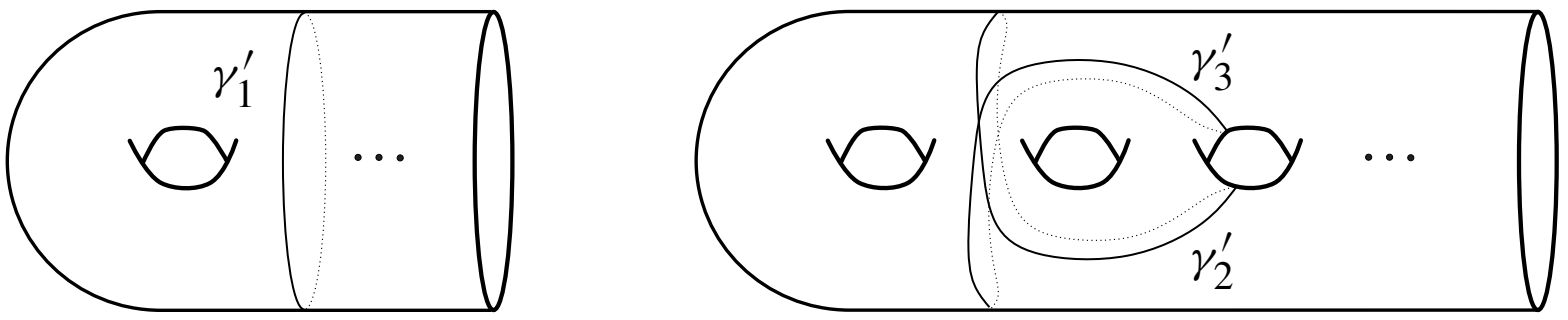

Figure 1: Left: the boundary curve $\gamma_1'$ of a genus-one subsurface with one boundary of the surface defining the Dehn twist $T_{\gamma_1'}$. Right: Simple closed curves $\gamma_2', \gamma_3'$ defining a homological genus-zero bounding pair map $B_0 := T_{\gamma_2'} T_{\gamma_3'}{}^{-1}$.

**Theorem D** *For $g \geq 6$, the order of the Euler class of the natural central extension*

$$0 \to \mathbb{Z} \to \mathrm{Ch}_{g,1} \to \mathrm{Ch}_{g,*} \to 1$$

*equals $\frac{1}{2}g(g-1)$ in $H^2(\mathrm{Ch}_{g,*};\mathbb{Z})$, and the abelianization of the Chillingworth subgroup $(\mathrm{Ch}_{g,*})^{ab} \cong H_1(\mathrm{Ch}_{g,*};\mathbb{Z})$ for the surface with a base point has a $\frac{1}{2}g(g-1)$-torsion element.*

## 2 Preliminaries

Let $\Sigma_{g,1}$ denote a connected, compact, oriented, genus-$g$ surface with one boundary. We choose a base point on the boundary of the surface $\Sigma_{g,1}$ and let $\{\alpha_1,\dots,\alpha_g,\beta_1,\dots,\beta_g\}$ be a free generating set of the fundamental group $\pi_1(\Sigma_{g,1})$ of the surface as shown in Figure 2.

Given two elements $\gamma_1$, $\gamma_2$ in the fundamental group of the surface $\pi = \pi_1(\Sigma_{g,1})$, their product $\gamma_1\gamma_2$ indicates that we traverse $\gamma_1$ first, then $\gamma_2$. The commutator $[\gamma_1,\gamma_2]$ is defined by $\gamma_1\gamma_2\gamma_1{}^{-1}\gamma_2{}^{-1}$.

Let $H = H_1(\Sigma_{g,1};\mathbb{Z})$ be the first integral homology group of the surface and $\cdot : H \otimes H \to \mathbb{Z}$ be the intersection form of the first homology of the surface. We choose a symplectic basis $\{a_1,\dots,a_g,b_1,\dots,b_g\}$ of $H$ as shown in Figure 3.

These elements are obtained through the Hurewicz homomorphism $\alpha_i \mapsto a_i$, $\beta_i \mapsto b_i$. The first integral cohomology group of the surface $H^* = H^1(\Sigma_{g,1};\mathbb{Z})$ is naturally isomorphic to the first homology group $H$ of the surface as $\mathrm{Sp}(2g,\mathbb{Z})$-modules by the Poincaré duality: $a_i \leftrightarrow b_i^*$, $b_i \leftrightarrow -a_i^*$. Using this, henceforth $H$ and $H^*$ will be freely identified. Let $\mathcal{M}_{g,1}$ be the mapping class group of the surface, which is defined as the isotopy classes of orientation-preserving self-diffeomorphisms of the surface that are pointwise identities on the boundary of the surface. That is,

$$\mathcal{M}_{g,1} := \mathrm{Diff}^{(+)}(\Sigma_{g,1}, \partial\Sigma_{g,1})/(\text{isotopies fixing the boundary pointwise}).$$

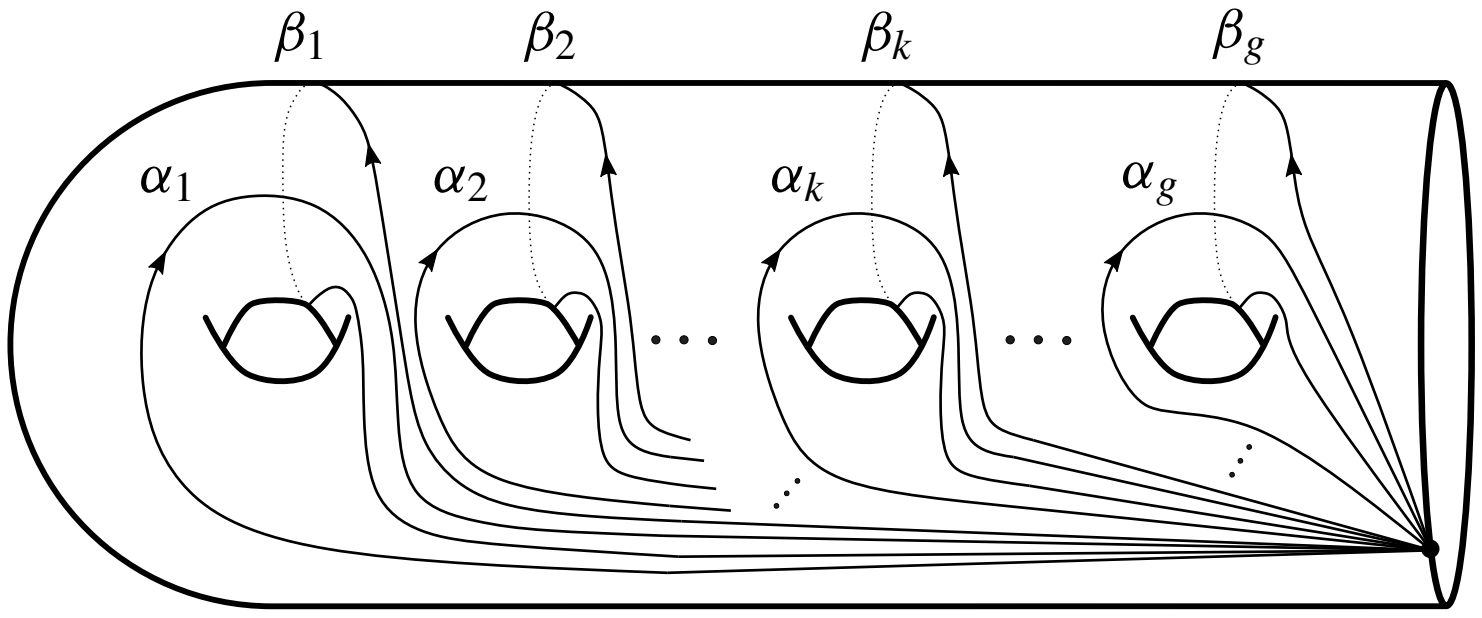

Figure 2: A generating system of the fundamental group of the surface $\pi_1(\Sigma_{g,1})$.

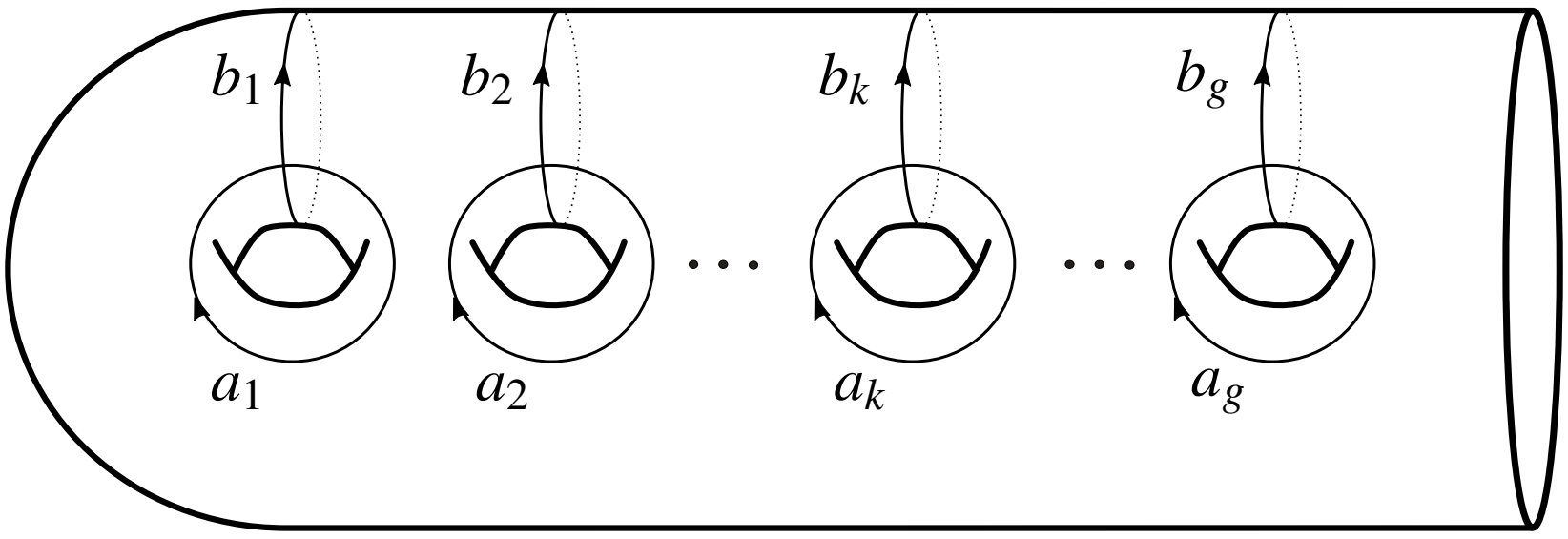

Figure 3: A symplectic basis of $H$.

A diffeomorphism that is the identity on the boundary is automatically orientation-preserving. The product $\varphi\psi$ in the mapping class group $\mathcal{M}_{g,1}$ indicates that we apply $\psi$ first, then $\varphi$. For a simple closed curve $C \subset \mathrm{Int}(\Sigma_{g,1})$, let $T_C$ be the (right-hand) Dehn twist along $C$.

## 2A Mapping class groups, fundamental groups, and Johnson homomorphisms

The action of the mapping class group on the fundamental group of the surface yields the Dehn–Nielsen representation $r : \mathcal{M}_{g,1} \to \mathrm{Aut}(\pi)$, which is known to be faithful. The mapping class group also acts naturally on the first integral homology group of the surface $H = H_1(\Sigma_{g,1}; \mathbb{Z})$ and this action preserves the intersection form of the surface. Hence, the mapping class acts on $H$ as the integral symplectic group $\mathrm{Sp}(H, \cdot) \cong \mathrm{Sp}(2g, \mathbb{Z})$ and this action $\rho : \mathcal{M}_{g,1} \to \mathrm{Sp}(2g, \mathbb{Z})$ is called the *symplectic representation*. It is known that the representation $\rho$ is surjective classically, and we summarize in the short exact sequence

$$1 \to \mathcal{I}_{g,1} \to \mathcal{M}_{g,1} \to \mathrm{Sp}(2g, \mathbb{Z}) \to 1,$$

where the kernel $\mathcal{I}_{g,1} := \mathrm{Ker}(\rho : \mathcal{M}_{g,1} \to \mathrm{Sp}(2g, \mathbb{Z}))$ of the symplectic representation is called the *Torelli group* of the mapping class group.

The Johnson homomorphism, initially defined by Johnson, provides an abelian quotient of the Torelli group and is equivariant under the action of the mapping class group (see Johnson [17; 18]). It has been developed by Morita and formalized as a graded Lie algebra homomorphism using the free Lie algebra generated by $H$ (see Morita [25; 29; 33]).

The mapping class group acts naturally on the nilpotent quotient of the fundamental group of the surface, denoted by $N_i := \pi/\Gamma_i$, where $\{\Gamma_i\}_{i\geq 1}$ is the lower central series of $\pi$, defined inductively by $\Gamma_1 := \pi$ and $\Gamma_{i+1} := [\Gamma_i, \pi]$.

**Definition 2.1** These actions on $\{N_i\}_{i\geq 1}$ define a filtration of the mapping class group, denoted by $\mathcal{M}_{g,1}[i] := \mathrm{Ker}(\mathcal{M}_{g,1} \to \mathrm{Aut}(N_i))$, called the Johnson filtration.

**Definition 2.2** The subgroup

$$\mathcal{K}_{g,1} := \langle \text{Dehn twists along bounding simple closed curves } (\textit{BSCC map}) \rangle$$

is called the Johnson kernel.

**Proposition 2.3** *We have*

$$\mathcal{M}_{g,1}[1] = \mathcal{M}_{g,1}, \quad \mathcal{M}_{g,1}[2] = \mathcal{I}_{g,1} = \mathrm{Ker}(\rho), \quad \mathcal{M}_{g,1}[3] = \mathcal{K}_{g,1}.$$

The last was shown by Johnson [20].

For $\varphi \in \mathcal{M}_{g,1}[i+1]$ and $\gamma \in \pi$, we have $\varphi(\gamma)\gamma^{-1} \in \Gamma_{i+1}$ by definition. Therefore, this defines a homomorphism $\mathcal{M}_{g,1}[i+1] \to \mathrm{Hom}(H, \Gamma_{i+1}/\Gamma_{i+2})$. The associated graded abelian group $\{\Gamma_i/\Gamma_{i+1}\}_{i\geq 1}$ of $\{\Gamma_i\}_{i\geq 1}$ admits a Lie algebra structure over $\mathbb{Z}$ via commutators on $\pi$. It is well known that the associated graded Lie algebra $\{\Gamma_i/\Gamma_{i+1}\}_{i\geq 1}$ is isomorphic to $\mathcal{L}_{g,1} = \{\mathcal{L}_{g,1}[i]\}_{i\geq 1}$, which is the free Lie algebra generated by $H$ over $\mathbb{Z}$, as a graded Lie algebra over $\mathbb{Z}$. For example, see [22]. Combining this with Poincaré duality, the homomorphism $\tau_{g,1}(i) : \mathcal{M}_{g,1}[i+1] \to H \otimes \mathcal{L}_{g,1}[i+1]$ is defined. By definition, we have $\mathrm{Ker}(\tau_{g,1}(i)) = \mathcal{M}_{g,1}[i+2]$.

Morita refined the target space using the structure of Lie algebra. Let

$$\mathfrak{h}_{g,1} = \{\mathfrak{h}_{g,1}(i)\}_{i\geq 1} := \left\{\mathrm{Ker}(H \otimes \mathcal{L}_{g,1}[i+1] \xrightarrow{\text{bracket}} \mathcal{L}_{g,1}[i+2])\right\}_{i\geq 1}$$

be the kernel of the bracket, which is a graded Lie subalgebra of $\mathrm{Hom}(H, \mathcal{L}_{g,1}) \cong \{H \otimes \mathcal{L}_{g,1}[i]\}_{i\geq 1}$.

**Theorem 2.4** (Morita [28; 29]) *The image* $\mathrm{Im}(\tau_{g,1}(i))$ *lies in* $\mathfrak{h}_{g,1}(i)$, *and* $\{\mathrm{Im}(\tau_{g,1}(i))\}_{i\geq 1}$ *is a graded Lie subalgebra of* $\mathfrak{h}_{g,1}$.

**Definition 2.5** (Morita) The homomorphism $\tau_{g,1}(i) : \mathcal{M}_{g,1}[i+1] \to \mathfrak{h}_{g,1}(i)$, known as the $i$-th Johnson homomorphism, is an $\mathcal{M}_{g,1}$-equivariant graded Lie algebra homomorphism

$$\{\tau_{g,1}(i)\}_{i\geq 1} : \{\mathcal{M}_{g,1}[i+1]/\mathcal{M}_{g,1}[i+2]\}_{i\geq 1} \to \{\mathfrak{h}_{g,1}(i)\}_{i\geq 1},$$

which is also called the Johnson homomorphism.

**Remark 2.6** Originally, Johnson defined it as $\tau_{g,1}(1) : \mathcal{I}_{g,1} \to \bigwedge^3 H$, where

$$\textstyle\bigwedge^3 H = \{x \wedge y \wedge z := x \otimes (y \wedge z) + y \otimes (z \wedge x) + z \otimes (x \wedge y) \mid x, y, z \in H\} \subset H \otimes \bigwedge^2 H \cong H \otimes \mathcal{L}_{g,1}[2],$$

and showed in [17] its surjectivity.

The above argument gives the Johnson filtrations and the Johnson homomorphisms for the mapping class group $\mathcal{M}_{g,*}$ of the surface with a base point and the mapping class group $\mathcal{M}_g$ of the closed surface without a base point.

**Definition 2.7** The homomorphism $\tau_{g,*}(i) : \mathcal{M}_{g,*}[i+1] \to \mathfrak{h}_{g,*}(i)$ is called the $i$-th Johnson homomorphism for $\mathcal{M}_{g,*}[i+1]$, where the target space $\mathfrak{h}_{g,*}(i)$ is defined as

$$\mathfrak{h}_{g,*} = \{\mathfrak{h}_{g,*}(i)\}_{i\geq 1} := \left\{\mathrm{Ker}(H \otimes \mathcal{L}_g[i+1] \xrightarrow{\text{bracket}} \mathcal{L}_g[i+2])\right\}_{i\geq 1},$$

where $\mathcal{L}_g := \mathcal{L}_{g,1}/\big(\omega_0 := \sum_{i=1}^{g}[a_i, b_i]\big)$. Similarly, the homomorphism $\tau_g(i) : \mathcal{M}_g[i+1] \to \mathfrak{h}_g(i)$ is called the $i$-th Johnson homomorphism for $\mathcal{M}_g[i+1]$, where the target space $\mathfrak{h}_g(i)$ is defined as

$$\mathfrak{h}_g := \mathfrak{h}_{g,*}/\mathcal{L}_g.$$

**Remark 2.8** By a result of Labute [22], the Lie algebra $\mathcal{L}_g$ is isomorphic to $\{\Gamma_i\pi_1(\Sigma_g)/\,\Gamma_{i+1}\pi_1(\Sigma_g)\}_{i\geq 1}$, where $\Gamma_i\pi_1(\Sigma)$ is the $i$-th term of the lower central series of $\pi_1(\Sigma_g)$.

**Remark 2.9** Originally, Johnson defined the first Johnson homomorphisms of these cases as $\tau_{g,*}(1) : \mathcal{I}_{g,*} \to \bigwedge^3 H$ and

$$\tau_g(1) : \mathcal{I}_g \to \textstyle\bigwedge^3 H/H := \bigwedge^3 H/\operatorname{Im}\Big(u : H \hookrightarrow \bigwedge^3 H,\ u(x) = \displaystyle\sum_{i=1}^{g} a_i \wedge b_i \wedge x\Big).$$

In this paper, for the sake of convenience, calculations using the first Johnson homomorphism are primarily performed using the original notation.

These Johnson homomorphisms commute with natural homomorphisms $\mathcal{M}_{g,1} \to \mathcal{M}_{g,*}$ induced by collapsing the boundary and $\mathcal{M}_{g,*} \to \mathcal{M}_g$ induced by forgetting the base point. There exists the commutative diagram

$$\begin{array}{ccccccccc}
1 & \longrightarrow & \mathcal{M}_{g,1}[i+2] & \longrightarrow & \mathcal{M}_{g,1}[i+1] & \xrightarrow{\tau_{g,1}(i)} & \mathfrak{h}_{g,1}(i) & \longrightarrow & 1 \\
 & & \downarrow & & \downarrow & & \downarrow & & \\
1 & \longrightarrow & \mathcal{M}_{g,*}[i+2] & \longrightarrow & \mathcal{M}_{g,*}[i+1] & \xrightarrow{\tau_{g,*}(i)} & \mathfrak{h}_{g,*}(i) & \longrightarrow & 1 \\
 & & \downarrow & & \downarrow & & \downarrow & & \\
1 & \longrightarrow & \mathcal{M}_{g}[i+2] & \longrightarrow & \mathcal{M}_{g}[i+1] & \xrightarrow{\tau_{g}(i)} & \mathfrak{h}_{g}(i) & \longrightarrow & 1
\end{array}$$

that commutes with the action of the mapping class group.

The short exact sequences

$$0 \to \mathbb{Z} \to \mathcal{M}_{g,1} \to \mathcal{M}_{g,*} \to 1,$$
$$0 \to \mathbb{Z} \to \mathcal{I}_{g,1} \to \mathcal{I}_{g,*} \to 1,$$
$$0 \to \mathbb{Z} \to \mathcal{K}_{g,1} \to \mathcal{K}_{g,*} \to 1$$

are induced by natural homomorphisms $\mathcal{M}_{g,1} \to \mathcal{M}_{g,*}$ and $\mathcal{M}_{g,*} \to \mathcal{M}_g$, where the second homomorphism from the left for each exact sequence is defined by $1 \mapsto T_\zeta$ and $\zeta$ is the boundary parallel loop of $\Sigma_{g,1}$. We also have the short exact sequences

$$1 \to \pi_1(\Sigma_g) \to \mathcal{M}_{g,*} \to \mathcal{M}_g \to 1,$$
$$1 \to \pi_1(\Sigma_g) \to \mathcal{I}_{g,*} \to \mathcal{I}_g \to 1,$$
$$1 \to [\pi_1(\Sigma_g), \pi_1(\Sigma_g)] \to \mathcal{K}_{g,*} \to \mathcal{K}_g \to 1.$$

The second homomorphism from the left for each exact sequence is called the push map defined by dragging the base point of the fundamental group along the element of the fundamental group. More generally, Asada and Kaneko showed in [1] that $\pi_1(\Sigma_g) \cap \mathcal{M}_{g,*}[i+1] = \Gamma_i \pi_1(\Sigma_g)$ and

$$1 \to \Gamma_i \pi_1(\Sigma_g) \to \mathcal{M}_{g,*}[i+1] \to \mathcal{M}_g[i+1] \to 1,$$

where $\Gamma_i G$ is the $i$-th term of the lower central series of $G$.

## 2B Tree diagrams and infinitesimal Dehn–Nielsen representations

The *infinitesimal Dehn–Nielsen representation*, introduced by Massuyeau [23], is an infinitesimal version of the Dehn–Nielsen representation $r : \mathcal{M}_{g,1} \to \mathrm{Aut}(\pi)$. It is described using an action on a certain complete Lie algebra defined by $\pi$, rather than the action on $\pi$ itself.

The target space of the infinitesimal Dehn–Nielsen representation is represented by $H$-labeled trees (see [12; 23]) called *tree diagrams*.

**Definition 2.10** A tree diagram is a finite, connected, unitrivalent graph whose trivalent vertices have cyclic order, and univalent vertices are colored by an element of $H$. The trivalent vertices of a tree diagram are called *nodes*, univalent vertices are called *leaves*, and the number of nodes in a tree diagram is called its *degree*.

**Definition 2.11** [12; 23] Define $\mathcal{T}_d(H)$ as the free abelian group generated by degree-$d$ tree diagrams modulo the relations

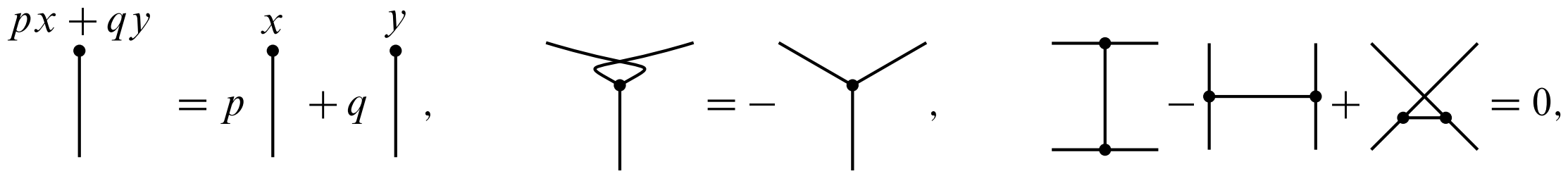

multilinearity AS relation (antisymmetry) IHX relation (Jacobi identity)

where $x, y \in H$ and $p, q \in \mathbb{Z}$. We define $\mathcal{T}(H) := \bigoplus_{d=1}^{\infty} \mathcal{T}_d(H)$, and $\widehat{\mathcal{T}(H)}$ as the degree completion of $\mathcal{T}(H)$. Similarly, we can define $\mathcal{T}(H_{\mathbb{Q}})$ over $\mathbb{Q}$ by taking the tensor product with $\mathbb{Q}$, giving us $\mathcal{T}(H_{\mathbb{Q}}) = \mathcal{T}(H) \otimes \mathbb{Q}$, and similarly for its completion $\widehat{\mathcal{T}(H_{\mathbb{Q}})}$, where the subscript $\mathbb{Q}$ means taking the tensor product $- \otimes \mathbb{Q}$.

Additionally, $\mathcal{T}(H)$ forms a graded Lie algebra over $\mathbb{Z}$, with the bracket $[\bullet, \bullet]_{\mathcal{T}}$ defined as

$$[P, Q]_{\mathcal{T}} := \sum_{\substack{v \in \mathrm{leaves}(P) \\ w \in \mathrm{leaves}(Q)}} (\mathrm{col}(P_v) \cdot \mathrm{col}(Q_w))(\text{graph obtained by gluing } P \text{ and } Q \text{ at } v \text{ and } w),$$

where $\mathrm{leaves}(P)$ is the set of leaves of $P$, $\mathrm{col}(P_v)$ is the color of the univalent vertex $v$, and $P_v$ is the rooted tree obtained by viewing $P$ as a rooted tree with root at vertex $v$. This bracket on $\mathcal{T}(H)$ is uniquely extended to the continuous bracket $[\bullet, \bullet]_{\widehat{\mathcal{T}}}$ on $\widehat{\mathcal{T}(H)}$. Then $(\widehat{\mathcal{T}(H)}, [\bullet, \bullet]_{\widehat{\mathcal{T}}})$ forms a complete graded Lie algebra over $\mathbb{Z}$. We can define similarly $(\widehat{\mathcal{T}(H_{\mathbb{Q}})}, [\bullet, \bullet]_{\widehat{\mathcal{T}}})$ over $\mathbb{Q}$.

The direct sum $\mathfrak{h}_{g,1} = \bigoplus_{i=1}^{\infty} \mathfrak{h}_{g,1}(i)$ of the target spaces of the Johnson homomorphisms forms a Lie subalgebra of $\bigoplus_{i=1}^{\infty} H \otimes \mathcal{L}_{g,1}[i+1]$, and there exists a Lie algebra homomorphism $\eta : \mathcal{T}(H) \to \mathfrak{h}_{g,1}$.

**Definition 2.12** The Lie algebra homomorphism $\eta : \mathcal{T}(H) \to \mathfrak{h}_{g,1}$ is defined by

$$\eta(P) := \sum_{v \in \mathrm{leaves}(P)} \mathrm{col}(P_v) \otimes \mathrm{brack}(P_v),$$

where $\mathrm{brack}(P_v)$ is the bracketification map defined by taking the iterated bracket as

brack ( 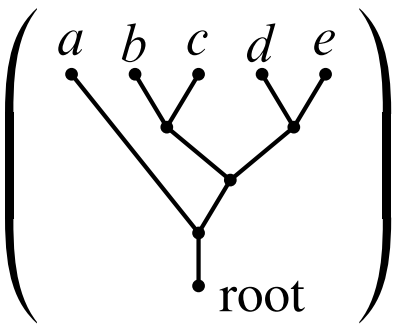

) $= \big[a, [[b, c], [d, e]]\big].$

For example, in the case of $d = 2$,

$\eta$ ( 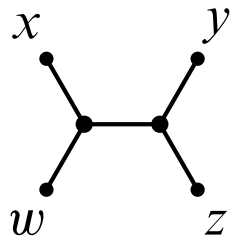

) $= x \otimes [[y, z], w] + y \otimes [z, [w, x]] + z \otimes [[w, x], y] + w \otimes [x, [y, z]].$

Especially, in the case of $d = 1$, we have the correspondence

$$\begin{matrix} x \wedge y \wedge z \\ \in \bigwedge^3 H \end{matrix} \longleftrightarrow \begin{matrix} x \otimes [y, z] + y \otimes [z, x] + z \otimes [x, y] \\ \in H \otimes \mathcal{L}_{g,1}[2] \end{matrix} \longleftrightarrow$$

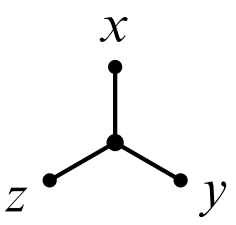

$\in \mathcal{T}_1(H).$

The homomorphism $\eta$ is not an isomorphism over $\mathbb{Z}$ (see [12]). However, if we take the tensor product with $\mathbb{Q}$, we have $\eta_{\mathbb{Q}} := \eta \otimes \mathrm{id}_{\mathbb{Q}}$, which is an isomorphism of Lie algebras over $\mathbb{Q}$ (see [14, Theorem 2]).

**Definition 2.13** (Massuyeau [23]) Let $\widehat{\mathcal{L}_{g,1\mathbb{Q}}}$ be the degree completion of $\mathcal{L}_{g,1\mathbb{Q}} = \bigoplus_{i=1}^{\infty} \mathcal{L}_{g,1\mathbb{Q}}[i]$. The *symplectic expansion* (*logansion*), introduced by Massuyeau, is a map $\theta : \pi \to \widehat{\mathcal{L}_{g,1\mathbb{Q}}}$ that satisfies the following conditions:

(1) $\theta : \pi \to (\widehat{\mathcal{L}_{g,1\mathbb{Q}}}, \star)$ is a group homomorphism, where $\star$ is the product defined by the *Baker–Campbell–Hausdorff series* (*BCH product*) with respect to the bracket $[\bullet, \bullet]$ of $\mathcal{L}_{g,1\mathbb{Q}}$.

(2) $\theta(x) = [x] + (\text{degree} \geq 2)$, where $[x] \in H$ is the image of $x \in \pi$ under the Hurewicz homomorphism.

(3) $\theta(\zeta) = -\omega := -\sum_{i=1}^{g} [a_i, b_i]$, where $\zeta \in \pi$ is the boundary parallel loop of $\Sigma_{g,1}$.

**Theorem 2.14** (Massuyeau [23, Lemma 2.16]) *There does exist a symplectic expansion* $\theta : \pi \to \widehat{\mathcal{L}_{g,1\mathbb{Q}}}$.

Before introducing the definition of the infinitesimal Dehn–Nielsen representation, we will introduce several maps. Let us fix a symplectic expansion $\theta : \pi \to \widehat{\mathcal{L}_{g,1\mathbb{Q}}}$.

First, we introduce the map $\varrho^{\theta} : \mathcal{I}_{g,1} \to \mathrm{IAut}_{\omega}(\widehat{\mathcal{L}_{g,1\mathbb{Q}}})$, where $\mathrm{IAut}_{\omega}(\widehat{\mathcal{L}_{g,1\mathbb{Q}}})$ is the automorphism group of $\widehat{\mathcal{L}_{g,1\mathbb{Q}}}$ whose elements induce the identity map on

$$\mathrm{gr}\mathcal{L}_{g,1\mathbb{Q}} = \bigoplus_{i=1}^{\infty} (\mathcal{L}_{g,1\mathbb{Q}}[i] / \mathcal{L}_{g,1\mathbb{Q}}[i+1])$$

and preserve the element $\omega = \sum_{i=1}^{g}[a_i, b_i] \in \mathcal{L}_{g,1}[2]$. The value $\varrho^\theta(f)$ for $f \in \mathcal{M}_{g,1}$ is induced by $x \mapsto f_*\theta(x)$.

Next, define $\log : \mathrm{IAut}_\omega(\widehat{\mathcal{L}_{g,1}}_{\mathbb{Q}}) \to (\mathrm{Der}^+_\omega(\widehat{\mathcal{L}_{g,1}}_{\mathbb{Q}}), \star)$, where the target space $\mathrm{Der}^+_\omega(\widehat{\mathcal{L}_{g,1}}_{\mathbb{Q}})$ is the space of derivations of $\widehat{\mathcal{L}_{g,1}}_{\mathbb{Q}}$ that strictly increase degrees and are trivial on $\omega$. The target space has a natural Lie algebra structure. If we define the group structure by the BCH product induced from this bracket, then the map log defined by the series induces a group isomorphism (the inverse map is exp defined by the series).

Finally, the map $\mathrm{Der}^+_\omega(\widehat{\mathcal{L}_{g,1}}_{\mathbb{Q}}) \to \widehat{\mathfrak{h}_{g,1}}$ is a canonical isomorphism induced by Poincaré duality. Specifically, it is the map induced by $D \mapsto \sum_{i=1}^{g}(b_i \otimes D(a_i) - a_i \otimes D(b_i))$.

Combining the above three maps with the map $\widehat{\eta_{\mathbb{Q}}^{-1}} : \widehat{\mathfrak{h}_{g,1}} \to \widehat{\mathcal{T}(H_{\mathbb{Q}})}$, we obtain the infinitesimal Dehn–Nielsen representation.

**Definition 2.15** The infinitesimal Dehn–Nielsen representation $r^\theta : \mathcal{I}_{g,1} \to (\widehat{\mathcal{T}(H_{\mathbb{Q}})}, \star)$ is defined as the composition of the homomorphisms

$$\mathcal{I}_{g,1} \xrightarrow{\varrho^\theta} \mathrm{IAut}_\omega(\widehat{\mathcal{L}_{g,1}}_{\mathbb{Q}}) \xrightarrow{\log} (\mathrm{Der}^+_\omega(\widehat{\mathcal{L}_{g,1}}_{\mathbb{Q}}), \star) \longrightarrow (\widehat{\mathfrak{h}_{g,1}}, \star) \longrightarrow (\widehat{\mathcal{T}(H_{\mathbb{Q}})}, \star), \qquad \mathcal{I}_{g,1} \overset{r^\theta}{\dashrightarrow} (\widehat{\mathcal{T}(H_{\mathbb{Q}})}, \star)$$

We also define its degree-$d$ part by composing with the projection $r^\theta_d : \mathcal{I}_{g,1} \xrightarrow{r^\theta} \widehat{\mathcal{T}(H_{\mathbb{Q}})} \to \mathcal{T}_d(H_{\mathbb{Q}})$; in particular, $\eta_{\mathbb{Q}} \circ r^\theta_i|_{\mathcal{M}_{g,1}[i+1]} : \mathcal{M}_{g,1}[i+1] \to \mathfrak{h}_{g,1\mathbb{Q}}(i)$ is nothing but the $i$-th Johnson homomorphism $\tau_{g,1}(i) : \mathcal{M}_{g,1}[i+1] \to \mathfrak{h}_{g,1}(i)$. Hence, the infinitesimal Dehn–Nielsen representation $r^\theta : \mathcal{I}_{g,1} \to \widehat{\mathcal{T}(H_{\mathbb{Q}})}$ on the $(i+1)$-st depth of the Johnson filtration $\mathcal{M}_{g,1}[i+1]$ is trivial up to degree-$(i-1)$-st part.

## 3 The action on the sets of homotopy classes of vector fields and the Chillingworth subgroups

Let $X$ be a nonsingular vector field on the surface $\Sigma_{g,1}$ and $\Xi(\Sigma_{g,1})$ be the set of homotopy classes of nonsingular vector fields on the surface. A homotopy class of nonsingular vector fields on a surface induces a trivialization of the unit tangent bundle $\mathrm{UT}\Sigma_{g,1} \xrightarrow{\cong} \Sigma_{g,1} \times S^1$ of the surface up to homotopy.

Let $\gamma$ be an oriented regular closed curve on the surface. The winding number of $\gamma$ with respect to $X$ denoted by $\omega_X(\gamma)$ is defined by the number of times its tangent transversely intersects with the section of the unit tangent bundle $\mathrm{UT}\Sigma_{g,1} \to \Sigma_{g,1}$ induced by $X$. Alternatively, we can compute the winding number by counting the points where the velocity vector is tangent to the vector field $X$, with the sign as shown in Figure 4.

The winding number function $\omega_X$ can be regarded as an element of $H^1(\mathrm{UT}\Sigma_{g,1}; \mathbb{Z})$. This element is characterized by the preimage of $1 \in H^1(S^1; \mathbb{Z})$ under the map $H^1(\mathrm{UT}\Sigma_{g,1}; \mathbb{Z}) \to H^1(S^1; \mathbb{Z})$. Conversely, for an arbitrary element $\omega \in H^1(\mathrm{UT}\Sigma_{g,1}; \mathbb{Z})$ which satisfies the condition, there exists a nonsingular vector field $X \in \Xi(\Sigma_{g,1})$ such that $\omega = \omega_X \in H^1(\mathrm{UT}\Sigma_{g,1}; \mathbb{Z})$; one can construct such an $X$ by considering $\Sigma_{g,1}$ as a disk with $2g$ attached 1-handles and specifying the vector field on each 1-handle. This correspondence $\Xi(\Sigma_{g,1}) \leftrightarrow \{\text{preimage of } 1\} \subset H^1(\mathrm{UT}\Sigma_{g,1};\mathbb{Z})$ is one-to-one.

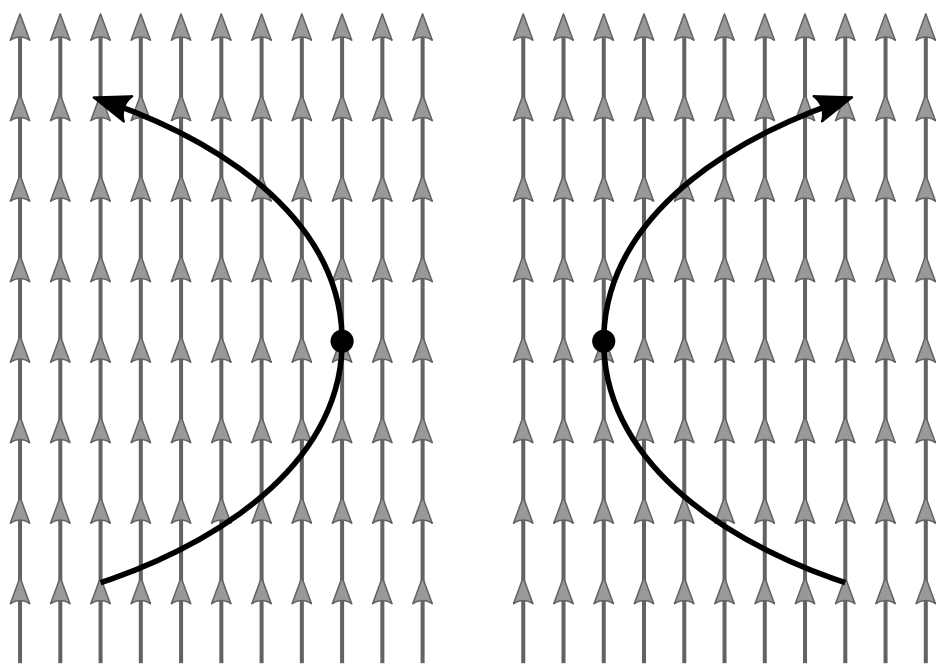

Figure 4: Signs of points where the velocity vector is tangent to the vector field.

The action of the mapping class group $\mathcal{M}_{g,1}$ of the surface on $\Xi(\Sigma_{g,1})$ is described by the $H^1(\Sigma_{g,1};\mathbb{Z})$-affine space structure and *the Chillingworth homomorphism*, which is defined using the winding number function $\omega_X$.

Let us fix a nonsingular vector field $X \in \Xi(\Sigma_{g,1})$. We recall the short exact sequence of the first cohomology

$$0 \to H^1(\Sigma_{g,1};\mathbb{Z}) \to H^1(\mathrm{UT}\Sigma_{g,1};\mathbb{Z}) \to H^1(S^1;\mathbb{Z}) \to 0,$$

which is equivariant under the action of the mapping class group.

**Definition 3.1** For a nonsingular vector field $X$, the Chillingworth homomorphism $e_X : \mathcal{M}_{g,1} \to H^1(\Sigma_{g,1};\mathbb{Z})$ is defined by the equality $e_X(f)([\gamma]) := \omega_X(f \circ \gamma) - \omega_X(\gamma)$.

The Chillingworth *homomorphism* is not a homomorphism but a crossed homomorphism. It satisfies $e_X(fg) = e_X(f) + (f^{-1})^* e_X(g)$. The kernel of the Chillingworth homomorphism $\mathrm{Ker}(e_X) := {e_X}^{-1}(0)$ is a subgroup of the mapping class group that preserve the chosen vector field $X$ up to homotopy. In particular, the Chillingworth homomorphism $e_X$ depends on the choice of a vector field $X$.

The construction of a crossed homomorphism on the mapping class group using the unit tangent bundle and the winding number was also proposed by Mikio Furuta. For details, see Morita [32]. Earle independently introduced an essentially identical crossed homomorphism using a different approach. For details, see [10].

Let us consider the restriction of the Chillingworth homomorphism to the Torelli subgroup. The restricted Chillingworth homomorphism $e_X|_{\mathcal{I}_{g,1}}$ is a homomorphism in the usual sense. Moreover, the restricted Chillingworth homomorphism does not depend on the choice of a nonsingular vector field on the surface. This is because, due to the short exact sequence, the difference $\omega_X - \omega_{X'}$ for different nonsingular vector fields $X$ and $X'$ can be expressed by an element $h \in H^1(\Sigma_{g,1};\mathbb{Z})$. From this, we have $e_X(f) - e_{X'}(f) = (f^{-1})^* h - h$, whose right-hand side is always zero on the Torelli group $\mathcal{I}_{g,1}$.

**Definition 3.2** The Chillingworth subgroup $\mathrm{Ch}_{g,1}$ is defined by the kernel $\mathrm{Ker}(e_X|_{\mathcal{I}_{g,1}})$ of the restricted Chillingworth homomorphism, and the Chillingworth subgroup of the surface with the base point $\mathrm{Ch}_{g,*}$ is defined similarly. We define the Chillingworth subgroup of the closed surface without a base point $\mathrm{Ch}_g$ as the image of the Chillingworth subgroup under the natural homomorphism $\mathcal{M}_{g,1} \to \mathcal{M}_{g,*} \to \mathcal{M}_g$.

Morita proved that $H^1(\mathcal{M}_{g,1}; H^*) \cong H^1(\mathcal{M}_{g,1}; H)$ is isomorphic to the infinite cyclic group $\mathbb{Z}$ in [26, Proposition 6.4] and that the crossed homomorphism (twisted 1-cocycle) $e_X$ is a generator of $H^1(\mathcal{M}_{g,1}; H^*)$ in [32, Proposition 4.1]. Hence, the Chillingworth subgroup is characterized as below.

**Proposition 3.3** (see [5; 7; 8; 41]) *The Chillingworth subgroup* $\mathrm{Ch}_{g,1}$ *has the following characterizations*:

(1) *the subgroup of the mapping class group whose elements preserve all nonsingular vector fields up to homotopy*;

(2) *the kernel* $\mathrm{Ker}(\mathcal{M}_{g,1} \curvearrowright \Xi(\Sigma_{g,1}))$ *of the action on the set of homotopy classes of nonsingular vector fields on the surface*;

(3) *the intersection of the kernel of a nontrivial crossed homomorphism with values in* $H$ *or* $H^*$ *and the Torelli group* $\mathcal{I}_{g,1}$;

(4) *the kernel* $\mathrm{Ker}(\mathcal{M}_{g,1} \curvearrowright H^1(\mathrm{UT}\Sigma_{g,1}; \mathbb{Z}))$ *of the action on the first cohomology of the unit tangent bundle of the surface*;

(5) *the kernel* $\mathrm{Ker}(\mathcal{M}_{g,1} \curvearrowright H_1(\mathrm{UT}\Sigma_{g,1}; \mathbb{Z}))$ *of the action on the first homology of the unit tangent bundle of the surface*;

(6) *the kernel* $\mathrm{Ker}(\mathcal{M}_{g,1} \curvearrowright \mathcal{H})$ *of the action on the Heisenberg group of the surface, where* $\mathcal{H}$ *is the Heisenberg group of the surface defined by* $\mathcal{H} = \mathbb{Z} \times H$ *as a set with the product defined by* $(n,x)(m,y) = (n+m+x\cdot y, x+y)$;

(7) *the kernel* $\mathrm{Ker}(\Phi_X : \mathcal{M}_{g,1} \to \mathrm{GL}(2g+1,\mathbb{Z}))$ *of Trapp's representation, which is defined by* $\Phi_X(f) = \begin{bmatrix} 1 & e_X(f) \\ 0 & \rho(f) \end{bmatrix}$.

First, there exists the following relationship among the Chillingworth subgroup, the Torelli group, and the Johnson kernel.

**Lemma 3.4** *For* $g \geq 3$, *we have* $\mathcal{K}_{g,1} \subsetneq \mathrm{Ch}_{g,1} \subsetneq \mathcal{I}_{g,1}$, *and for* $g = 2$, *we have* $\mathcal{K}_{2,1} = \mathrm{Ch}_{2,1} \subsetneq \mathcal{I}_{2,1}$. *Similarly, for* $g \geq 3$, *we have* $\mathcal{K}_{g,*} \subsetneq \mathrm{Ch}_{g,*} \subsetneq \mathcal{I}_{g,*}$ *and* $\mathcal{K}_g \subsetneq \mathrm{Ch}_g \subsetneq \mathcal{I}_g$, *and for* $g = 2$, *we have* $\mathcal{K}_{2,*} = \mathrm{Ch}_{2,*} \subsetneq \mathcal{I}_{g,*}$ *and* $\mathcal{K}_2 = \mathrm{Ch}_2 = \mathcal{I}_2$.

Furthermore, $\mathrm{Ch}_g$ is a finite-index normal subgroup of $\mathcal{I}_g$.

**Lemma 3.5** *For* $g \geq 3$, $\mathrm{Ch}_g$ *is a normal subgroup of index* $(g-1)^{2g}$ *in* $\mathcal{I}_g$, *and the quotient* $\mathcal{I}_g/\mathrm{Ch}_g$ *is isomorphic to* $(\mathbb{Z}/(g-1)\mathbb{Z})^{2g}$.

Additionally, the relationships between the Chillingworth subgroups in each case can be summarized by the following short exact sequences.

**Proposition 3.6** *There exist two short exact sequences*

$$0 \to \mathbb{Z} \to \mathrm{Ch}_{g,1} \to \mathrm{Ch}_{g,*} \to 1, \quad 1 \to [\pi_1(\Sigma_g), \pi_1(\Sigma_g)] \to \mathrm{Ch}_{g,*} \to \mathrm{Ch}_g \to 1.$$

Lemmas 3.4–3.5 and Proposition 3.6 can be seen from the relationship between the Chillingworth subgroup and the Johnson homomorphism, which will be explained in the next subsection.

## 3A The first Johnson homomorphisms and the Chillingworth subgroups

Johnson [17] introduced the element $t_f \in H$ as the Poincaré dual of the value of the Chillingworth homomorphism. It is characterized by the property that $x \cdot t_f = e_X(f)(x)$ for all $x \in H_1(\Sigma_{g,1}; \mathbb{Z})$. Here, $t_\bullet$ is called the *Chillingworth class*. Johnson proved that the Chillingworth class factors through the first Johnson homomorphism.

**Lemma 3.7** (Johnson [17, Theorem 2]) *The diagram*

$$\begin{array}{ccc} \mathcal{I}_{g,1} & \xrightarrow{\tau_{g,1}(1)} & \bigwedge^3 H \\ & {\scriptstyle t}\searrow & \downarrow {\scriptstyle 2C_3} \\ & & H \end{array}$$

*is $\mathcal{M}_{g,1}$-equivariant and commutative. Here, the $\mathrm{Sp}(2g,\mathbb{Z})$-equivariant homomorphism $C_3 : \bigwedge^3 H \to H$ is defined by $x \wedge y \wedge z \mapsto (x \cdot y)z + (y \cdot z)x + (z \cdot x)y$ and called the* ***contraction****.*

Moreover, Morita [30] constructed an extension of the Johnson homomorphism as a crossed homomorphism to the mapping class group $\mathcal{M}_{g,1} \to \frac{1}{2} \bigwedge^3 H$. By composing this extension with $2C_3$, one can also obtain a crossed homomorphism on $\mathcal{M}_{g,1}$.

Here, we introduce certain elements of a Torelli group that will appear in subsequent discussions.

**Definition 3.8** For two disjoint nonseparating simple closed curves $\gamma_1$ and $\gamma_2$ on the surface $\Sigma_{g,1}$, when there exists a subsurface with genus $h$ with the boundary components equal to $\gamma_1 \cup \gamma_2$, we call the map $\mathrm{BP}(\gamma_1, \gamma_2) := T_{\gamma_1} {T_{\gamma_2}}^{-1}$ the *genus-$h$ bounding pair map* (*BP map*), which is an element of the Torelli group.

To prove Lemmas 3.4–3.5 and Proposition 3.6, we introduce some calculation formulas.

**Proposition 3.9** (Johnson [17, Lemma 4A]) *Let $\{x_i, y_i\}_{i=1,\dots,h}$ be a symplectic basis of the first homology group of the subsurface defining a genus-$h$ BP map $\mathrm{BP}(\gamma_1, \gamma_2)$. Then, we have*

$$\tau_{g,1}(1)(\mathrm{BP}(\gamma_1, \gamma_2)) = \sum_{i=1}^{h} x_i \wedge y_i \wedge [\gamma_1],$$

*where $\gamma_1$ and $\gamma_2$ are endowed with the orientation induced from the subsurface.*

Furthermore, using $C_3(x_i \wedge y_i \wedge [\gamma_1]) = [\gamma_1]$, we obtain the following:

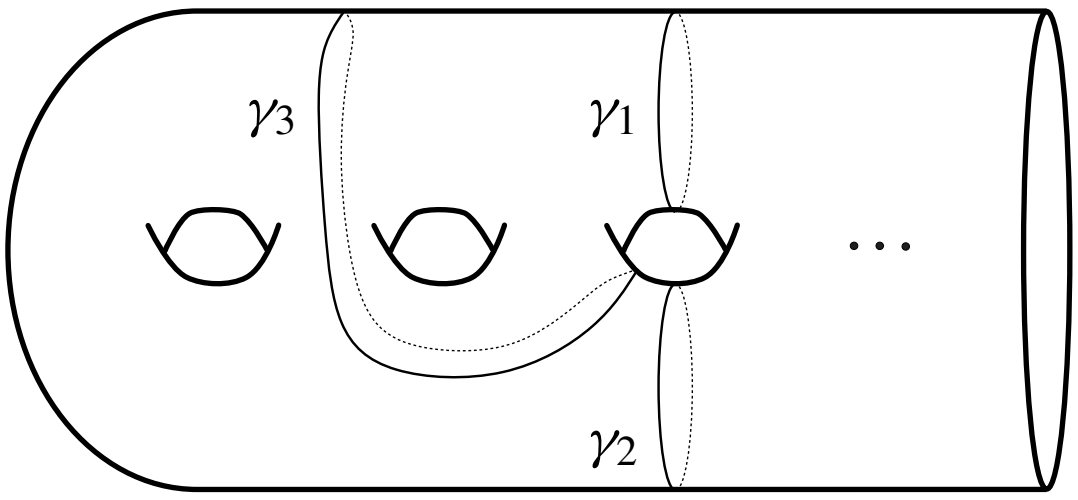

Figure 5: Some simple closed curves on the surface defining some BP maps.

**Lemma 3.10** *Let* $\mathrm{BP}(\gamma_1,\gamma_2)$ *be a genus-h BP map. We have*

$$t_{\mathrm{BP}(\gamma_1,\gamma_2)} = 2C_3 \circ \tau_{g,1}(1)(\mathrm{BP}(\gamma_1,\gamma_2)) = 2h[\gamma_1].$$

Now, we prove Lemma 3.4.

**Proof** For $g=2$, the contraction $C_3 : \bigwedge^3 H \to H$ is an isomorphism, and its kernel is trivial. Hence, the conditions $\tau_{2,1}(1)(f)=0$ and $t_f=0$ are equivalent, which implies that $\mathrm{Ch}_{2,1}=\mathcal{K}_{2,1}$. It follows that $\mathrm{Ch}_{2,*}=\mathcal{K}_{2,*}$ and $\mathrm{Ch}_2=\mathcal{K}_2$. Moreover, in the case of genus-two closed surfaces without a base point, the target space of the first Johnson homomorphism $\bigwedge^3 H/H$ is trivial. Therefore $\mathcal{K}_2=\mathcal{I}_2$; in particular, $\mathcal{K}_2=\mathrm{Ch}_2=\mathcal{I}_2$. Next, consider the case of $g\geq 3$ with one boundary component. Consider a genus-one BP map. This element is contained in $\mathcal{I}_{g,1}$ but not contained in $\mathrm{Ch}_{g,1}$, as its value under $t=2C_3\circ\tau_{g,1}(1)$ is nontrivial, as shown in Lemma 3.10. Therefore, we have $\mathrm{Ch}_{g,1}\subsetneq\mathcal{I}_{g,1}$. For $g\geq 3$, let us consider the element $\mathrm{BP}(\gamma_1,\gamma_2)\,\mathrm{BP}(\gamma_1,\gamma_3)^{-2}$ as in Figure 5. This element is contained in $\mathrm{Ch}_{g,1}$ but not in $\mathcal{K}_{g,1}$. Specifically, $2C_3\circ\tau_{g,1}(1)(\mathrm{BP}(\gamma_1,\gamma_2)\,\mathrm{BP}(\gamma_1,\gamma_3)^{-2})=2C_3(a_1\wedge b_1\wedge b_3-a_2\wedge b_2\wedge b_3)=0$, indicating this element is contained in the Chillingworth subgroup $\mathrm{Ch}_{g,1}$. However, $\tau_{g,1}(1)(\mathrm{BP}(\gamma_1,\gamma_2)\,\mathrm{BP}(\gamma_1,\gamma_3)^{-2})=a_1\wedge b_1\wedge b_3-a_2\wedge b_2\wedge b_3$ is nonzero, implying this element is not contained in the Johnson kernel $\mathcal{K}_{g,1}$. Therefore, we have $\mathcal{K}_{g,1}\subsetneq\mathrm{Ch}_{g,1}$. The same argument can be applied for $\mathrm{Ch}_{g,*}$ for $g\geq 2$, and $\mathrm{Ch}_g$ for $g\geq 3$ cases. □

Before Lemma 3.5, we discuss Proposition 3.6. Johnson [17] discusses the kernel of the homomorphism $\mathcal{I}_{g,1}\to\mathcal{I}_g$. We obtain the following:

**Lemma 3.11** (Johnson [17]) *By composing the push map* $\pi_1(\Sigma_g)\hookrightarrow\mathcal{I}_{g,*}$ *with the first Johnson homomorphism* $\tau_{g,*}(1):\mathcal{I}_{g,*}\to\bigwedge^3 H$, *we obtain*

$$\tau_{g,*}(1)(\gamma) = -\sum_{i=1}^{g} a_i\wedge b_i\wedge[\gamma],$$

*for* $\gamma\in\pi_1(\Sigma_g)$. *In particular,* $t_\gamma=-2(g-1)[\gamma]$.

This shows that $t \bmod (2g-2):\mathcal{I}_g\to H\otimes(\mathbb{Z}/(2g-2)\mathbb{Z})$ is well defined.

Using this lemma, we will now proceed to prove Proposition 3.6.

**Proof** Since

$$[\pi_1(\Sigma_g),\pi_1(\Sigma_g)]=\mathcal{K}_{g,*}\cap\pi_1(\Sigma_g)\subset \mathrm{Ch}_{g,*}\cap\pi_1(\Sigma_g),$$

we only need to prove that

$$[\pi_1(\Sigma_g),\pi_1(\Sigma_g)]\supset \mathrm{Ch}_{g,*}\cap\pi_1(\Sigma_g).$$

Let $\gamma$ be an element of $\pi_1(\Sigma_g)\cap\mathrm{Ch}_*$. Since $\gamma$ is contained in $\mathrm{Ch}_{g,*}$, we have $t_\gamma=0$. From Lemma 3.11, we have $-2(g-1)[\gamma]=0$ in $H=H_1(\Sigma_{g,*};\mathbb{Z})\cong H_1(\Sigma_g;\mathbb{Z})\cong\pi_1(\Sigma_g)/[\pi_1(\Sigma_g),\pi_1(\Sigma_g)]$. Since $H$ is a free abelian group, $[\gamma]=0$ in $\pi_1(\Sigma_g)/[\pi_1(\Sigma_g),\pi_1(\Sigma_g)]$. Therefore $\gamma\in[\pi_1(\Sigma_g),\pi_1(\Sigma_g)]$. □

We denote the kernel of the contraction $\mathrm{Ker}(C_3)\subset\bigwedge^3 H$ as $U$. Note that $U$ is a rank-$\left(\binom{2g}{3}-2g\right)$ free abelian group and a $\mathrm{Sp}(2g,\mathbb{Z})$-submodule of $\bigwedge^3 H$. We denote the image of $U$ under the natural homomorphism $U\hookrightarrow\bigwedge^3 H\to\bigwedge^3 H/H$ as $\overline{U}$. By definition, these coincide with the images of the Chillingworth subgroups under the Johnson homomorphisms: $\tau_{g,1}(1)(\mathrm{Ch}_{g,1})=\tau_{g,*}(1)(\mathrm{Ch}_{g,*})=U$ and $\tau_g(1)(\mathrm{Ch}_g)=\overline{U}$.

Finally, we prove Lemma 3.5, which follows from Lemma 3.12. Before stating the lemma, we define the map $v: H\oplus U\to\bigwedge^3 H$ as

$$v: H\oplus U\to\bigwedge^3 H,\quad (x,Y)\mapsto\Big(\sum_{i=1}^{g} a_i\wedge b_i\wedge x\Big)+Y.$$

**Lemma 3.12** *For $g\ge 3$, the quotient*

$$\mathcal{I}_g/\mathrm{Ch}_g\cong\big(\textstyle\bigwedge^3 H/H\big)/\overline{U}=\mathrm{Coker}\big(v: H\oplus U\to\bigwedge^3 H\big)$$

*is isomorphic to $(\mathbb{Z}/(g-1)\mathbb{Z})^{2g}$.*

**Proof** Let us take a basis of $U$ as

(i) $a_i\wedge a_j\wedge a_k$, $b_i\wedge b_j\wedge b_k$ for distinct $i,j,k$,

(ii) $a_i\wedge a_j\wedge b_k$, $a_i\wedge b_j\wedge b_k$ for distinct $i,j,k$,

(iii) $a_1\wedge a_2\wedge b_2-a_1\wedge a_i\wedge b_i$ for $i\ge 3$, $a_j\wedge a_1\wedge b_1-a_j\wedge a_i\wedge b_i$ for $i\ge 3$, $j\ge 2$, $i\ne j$, $b_1\wedge a_2\wedge b_2-a_1\wedge b_i\wedge b_i$ for $i\ge 3$, and $b_j\wedge a_1\wedge b_1-b_j\wedge a_i\wedge b_i$ for $i\ge 3$, $j\ge 2$, $i\ne j$,

and take a basis of $\bigwedge^3 H$ as (i), (ii), (iii), and

(iv) $a_i\wedge a_j\wedge b_j$ for $i\ne j$, $b_i\wedge a_j\wedge b_j$ for $i\ne j$.

The representation matrix of $v: H\oplus U\to\bigwedge^3 H$ with respect to the above basis is

$$(I_{\binom{g}{3}})^{\oplus 2}\oplus(I_{g\binom{g-1}{2}})^{\oplus 2}\oplus\begin{bmatrix}1&1&1&1&\cdots&1\\1&-1&&&&0\\1&&-1&&&\vdots\\1&&&-1&&\\\vdots&&&&\ddots&0\\1&0&\cdots&&0&-1\end{bmatrix}^{\oplus 2g},$$

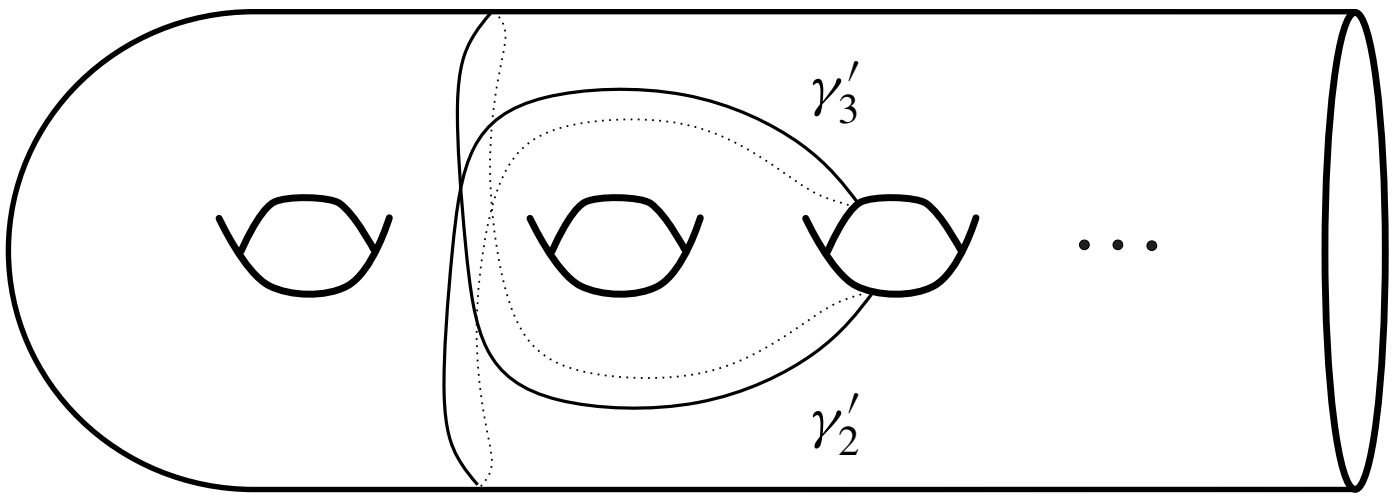

Figure 6: Some simple closed curves on the surface defining the element $B_0$.

where $I_n$ is the identity matrix of size $n\times n$ and the rightmost matrix is of $(g-1)\times(g-1)$ size. We compute the invariant factor of it as

$$(I_{\binom{g}{3}})^{\oplus 2} \oplus (I_{g\binom{g-1}{2}})^{\oplus 2} \oplus \begin{bmatrix} g-1 & & & \\ & 1 & & \\ & & \ddots & \\ & & & 1 \end{bmatrix}^{\oplus 2g}.$$

Hence, the quotient $\mathcal{I}_g/\mathrm{Ch}_g$ is isomorphic to $(\mathbb{Z}/(g-1)\mathbb{Z})^{2g}$. □

Lemma 3.12 is restated as follows.

**Proposition 3.13** *By Lemma 3.11, the map* $t \bmod (2g-2) : \mathcal{I}_g \to H \otimes \mathbb{Z}/(2g-2)\mathbb{Z}$ *is a well-defined homomorphism, with its kernel being* $\mathrm{Ch}_g$, *and its image being* $2H \otimes \mathbb{Z}/(2g-2)\mathbb{Z} \cong H \otimes \mathbb{Z}/(g-1)\mathbb{Z}$.

As stated previously, the composition $U \to \bigwedge^3 H \to \bigwedge^3 H/H$ is not an isomorphism. However, if we take the tensor product with $\mathbb{Q}$, then the composition

$$U \otimes \mathbb{Q} \to \left(\bigwedge^3 H\right) \otimes \mathbb{Q} \to \left(\bigwedge^3 H/H\right) \otimes \mathbb{Q}$$

becomes an isomorphism as $\mathrm{Sp}(2g,\mathbb{Q})$-modules. We use the notation $\bigwedge^3 H_{\mathbb{Q}} := \left(\bigwedge^3 H\right)\otimes\mathbb{Q} = \bigwedge^3(H_{\mathbb{Q}})$, $U_{\mathbb{Q}} := U \otimes \mathbb{Q}$ and so forth.

**Proposition 3.14** *For* $g \geq 4$, *Chillingworth subgroups* $\mathrm{Ch}_{g,1}$, $\mathrm{Ch}_{g,*}$ *and* $\mathrm{Ch}_g$ *are normally generated by one element and the Johnson kernel in the full mapping class group.*

**Proof** We consider the exact sequence

$$1 \to \mathcal{K}_{g,1} \to \mathrm{Ch}_{g,1} \to U \to 1$$

induced by the Johnson homomorphism for the Chillingworth subgroup. The Chillingworth subgroup $\mathrm{Ch}_{g,1}$ is generated by $\mathcal{K}_{g,1}$ together with lifts of elements of $U$ under the surjective homomorphism $\tau_{g,1}(1) : \mathrm{Ch}_{g,1} \to U$. Let us take the conjugacy class of a certain element $B_0 := \mathrm{BP}(\gamma_2', \gamma_3') := T_{\gamma_2'} T_{\gamma_3'}{}^{-1}$ as shown in Figure 6 (which we call a *homological genus-zero* (*or one minus one*) *bounding pair map*). The image of this conjugacy class under the Johnson homomorphism is surjective onto $U$. Equivalently, $U$ is generated by $\tau_{g,1}(1)(B_0) = a_1 \wedge b_1 \wedge b_3 - a_2 \wedge b_2 \wedge b_3$ as an $\mathrm{Sp}(2g,\mathbb{Z})$-module. To show this, it suffices to construct all the basis elements of $U$ given in the proof of Lemma 3.12(i), (ii), and (iii) by applying appropriate elements of $\mathrm{Sp}(2g,\mathbb{Z})$ to $a_1\wedge b_1\wedge b_3 - a_2\wedge b_2\wedge b_3$. By suitably permuting the indices and applying the matrices determined by $a_i\mapsto b_i\mapsto -a_i$, we can construct all the elements of type (iii). Next, if we subtract the result of applying

$$\begin{cases} a_1\mapsto a_1+b_1-b_4,\\ a_4\mapsto a_4+b_4-b_1,\\ a_i\mapsto a_i \quad (i\neq 1,4),\\ b_i\mapsto b_i \end{cases}$$

to the original element $a_1\wedge b_1\wedge b_3 - a_2\wedge b_2\wedge b_3$, we obtain $b_4\wedge b_1\wedge b_3$. From this element, we can similarly construct elements of types (i) and (ii) using the same argument. Therefore, $\mathrm{Ch}_{g,1}$ is normally generated by $B_0$ and the Johnson kernel. For $\mathrm{Ch}_{g,*}$ and $\mathrm{Ch}_g$, we obtain similar results via the natural surjective homomorphisms $\mathrm{Ch}_{g,1}\to\mathrm{Ch}_{g,*}$ and $\mathrm{Ch}_{g,*}\to\mathrm{Ch}_g$. □

## 4 Proof of Theorem B

By the general theory of representation, a finite-dimensional polynomial representation of the rational symplectic group $\mathrm{Sp}(2g,\mathbb{Q})$ corresponds bijectively to those of $\mathrm{Sp}(2g,\mathbb{C})$ and the Lie algebra $\mathfrak{sp}(2g,\mathbb{C})$. These representations are parametrized by Young diagrams. We use a notation in conformity to Fulton–Harris [13].

We denote the one-dimensional trivial representation $\mathbb{Q}$ by $[0]_{\mathrm{Sp}}$, and the natural representation $H_{\mathbb{Q}}$ by $[1]_{\mathrm{Sp}}$. For a Young diagram corresponding to $n_1\geq n_2\geq\cdots\geq n_l\geq 1$, $l\leq g$, we define $[n_1n_2\cdots n_l]_{\mathrm{Sp}}$ as below:

- Let $m_1\geq m_2\geq\cdots\geq m_k$ be the transpose of $n_1\geq n_2\geq\cdots\geq n_l\geq l$.
- The vector

$$(a_1\wedge a_2\wedge\cdots\wedge a_{m_1})\otimes(a_1\wedge a_2\wedge\cdots\wedge a_{m_2})\otimes\cdots\otimes(a_1\wedge a_2\wedge\cdots\wedge a_{m_k})$$

  within $\left(\bigwedge^{m_1} H_{\mathbb{Q}}\right)\otimes\left(\bigwedge^{m_2} H_{\mathbb{Q}}\right)\otimes\cdots\otimes\left(\bigwedge^{m_k} H_{\mathbb{Q}}\right)$ generates an irreducible subrepresentation.
- This irreducible representation is denoted by $[n_1n_2\cdots n_l]_{\mathrm{Sp}}$, and the vector

$$(a_1\wedge a_2\wedge\cdots\wedge a_{m_1})\otimes(a_1\wedge a_2\wedge\cdots\wedge a_{m_2})\otimes\cdots\otimes(a_1\wedge a_2\wedge\cdots\wedge a_{m_k})$$

  is called the *highest weight vector* of $[n_1n_2\cdots n_l]_{\mathrm{Sp}}$.

We abbreviate $[2211]_{\mathrm{Sp}}$, $[111111]_{\mathrm{Sp}}$ and so forth as $[2^21^2]_{\mathrm{Sp}}$, $[1^6]_{\mathrm{Sp}}$ and so forth.

These representations are naturally isomorphic to their dual representation. For example, $H_{\mathbb{Q}}^*$ and its dual $H_{\mathbb{Q}}$ are isomorphic as representations of $\mathrm{Sp}(2g,\mathbb{Q})$ via the Poincaré duality, and are denoted by $[1]_{\mathrm{Sp}}$. We have $\bigwedge^3 H_{\mathbb{Q}} = [1^3]_{\mathrm{Sp}}\oplus[1]_{\mathrm{Sp}}$ and $\bigwedge^3 H_{\mathbb{Q}}/H_{\mathbb{Q}}\cong U_{\mathbb{Q}} = [1^3]_{\mathrm{Sp}}$.

The following proposition follows from the irreducibility of $U_{\mathbb{Q}}=[1^3]_{\mathrm{Sp}}$.

**Proposition 4.1** *For $g \geq 3$,*

$$(\tau_{g,1}(1))_* : H_1(\mathrm{Ch}_{g,1};\mathbb{Q}) \to H_1(U;\mathbb{Q})$$

*is surjective and*

$$(\tau_{g,1}(1))^* : H^1(U;\mathbb{Q}) \to H^1(\mathrm{Ch}_{g,1};\mathbb{Q})$$

*is injective. The same holds for the* $\mathrm{Ch}_{g,*}$ *and* $\mathrm{Ch}_g$ *cases.*

Hain studied the homomorphism $(\tau_g(1))^* : H^2\big(\bigwedge^3 H/H;\mathbb{Q}\big) \to H^2(\mathcal{I}_g;\mathbb{Q})$ between the second rational cohomology induced by the Johnson homomorphism and determined the kernel of this map as $\mathrm{Sp}(2g,\mathbb{Q})$-modules using representation theory.

**Lemma 4.2** (Hain [15, Lemma 10.2]) *For $g \geq 3$,*

$$\begin{aligned} H^2\big(\textstyle\bigwedge^3 H/H;\mathbb{Q}\big) &\cong H_2\big(\textstyle\bigwedge^3 H/H;\mathbb{Q}\big) \cong H^2(U;\mathbb{Q}) \cong H_2(U;\mathbb{Q}) \cong \textstyle\bigwedge^2 U_{\mathbb{Q}} \\ &= \begin{cases} [0]_{\mathrm{Sp}} \oplus [2^2]_{\mathrm{Sp}} \oplus [1^2]_{\mathrm{Sp}} \oplus [2^2 1^2]_{\mathrm{Sp}} \oplus [1^4]_{\mathrm{Sp}} \oplus [1^6]_{\mathrm{Sp}} & (g \geq 6), \\ [0]_{\mathrm{Sp}} \oplus [2^2]_{\mathrm{Sp}} \oplus [1^2]_{\mathrm{Sp}} \oplus [2^2 1^2]_{\mathrm{Sp}} \oplus [1^4]_{\mathrm{Sp}} & (g = 5), \\ [0]_{\mathrm{Sp}} \oplus [2^2]_{\mathrm{Sp}} \oplus [1^2]_{\mathrm{Sp}} \oplus [2^2 1^2]_{\mathrm{Sp}} & (g = 4), \\ [0]_{\mathrm{Sp}} \oplus [2^2]_{\mathrm{Sp}} & (g = 3) \end{cases} \end{aligned}$$

*holds as* $\mathrm{Sp}(2g,\mathbb{Q})$*-modules.*

**Theorem 4.3** (Hain [15]) *For $g \geq 3$,*

$$\mathrm{Ker}\big((\tau_g(1))^* : H^2\big(\textstyle\bigwedge^3 H/H;\mathbb{Q}\big) \to H^2(\mathcal{I}_g;\mathbb{Q})\big) = [0]_{\mathrm{Sp}} \oplus [2^2]_{\mathrm{Sp}}$$

*holds as* $\mathrm{Sp}(2g,\mathbb{Q})$*-modules.*

Moreover, the dual of the preceding theorem implies that the image of the homomorphism

$$(\tau_g(1))_* : H_2(\mathcal{I}_g;\mathbb{Q}) \to H_2\big(\textstyle\bigwedge^3 H/H;\mathbb{Q}\big)$$

between the second rational homology induced by the Johnson homomorphism is decomposed as $\mathrm{Sp}(2g,\mathbb{Q})$-modules as follows:

**Theorem 4.4** (Hain [15]) *For $g \geq 3$,*

$$\mathrm{Im}\big((\tau_g(1))_* : H_2(\mathcal{I}_g;\mathbb{Q}) \to H_2\big(\textstyle\bigwedge^3 H/H;\mathbb{Q}\big)\big) = \begin{cases} [1^2]_{\mathrm{Sp}} \oplus [2^2 1^2]_{\mathrm{Sp}} \oplus [1^4]_{\mathrm{Sp}} \oplus [1^6]_{\mathrm{Sp}} & (g \geq 6), \\ [1^2]_{\mathrm{Sp}} \oplus [2^2 1^2]_{\mathrm{Sp}} \oplus [1^4]_{\mathrm{Sp}} & (g = 5), \\ [1^2]_{\mathrm{Sp}} \oplus [2^2 1^2]_{\mathrm{Sp}} & (g = 4), \\ \{0\} & (g = 3) \end{cases}$$

*holds as* $\mathrm{Sp}(2g,\mathbb{Q})$*-modules.*

For $g \geq 3$, the homomorphism $\mathrm{Ch}_{g,1} \hookrightarrow \mathcal{I}_{g,1} \to \mathcal{I}_{g,*} \to \mathcal{I}_g$ induces the $\mathcal{M}_{g,1}$-equivariant commutative diagram

$$\begin{array}{ccccc}
\bigwedge^2 U_\mathbb{Q} & \xrightarrow{\cong} & H^2(\bigwedge^3 H/H;\mathbb{Q}) & \xrightarrow{(\tau_g(1))^*} & H^2(\mathcal{I}_g;\mathbb{Q}) \\
\big\downarrow & & \big\downarrow & & \big\downarrow \\
\bigwedge^2 H_\mathbb{Q} \oplus (H_\mathbb{Q} \otimes U_\mathbb{Q}) \oplus \bigwedge^2 U_\mathbb{Q} & \xrightarrow{\cong} & H^2(\bigwedge^3 H;\mathbb{Q}) & \xrightarrow{(\tau_{g,*}(1))^*} & H^2(\mathcal{I}_{g,*};\mathbb{Q}) \\
\big\| & & \big\downarrow & & \big\downarrow \\
\bigwedge^2 H_\mathbb{Q} \oplus (H_\mathbb{Q} \otimes U_\mathbb{Q}) \oplus \bigwedge^2 U_\mathbb{Q} & \xrightarrow{\cong} & H^2(\bigwedge^3 H;\mathbb{Q}) & \xrightarrow{(\tau_{g,1}(1))^*} & H^2(\mathcal{I}_{g,1};\mathbb{Q}) \\
\big\downarrow & & \big\downarrow & & \big\downarrow \\
\bigwedge^2 U_\mathbb{Q} & \xrightarrow{\cong} & H^2(U;\mathbb{Q}) & \xrightarrow{(\tau_{g,1}(1))^*} & H^2(\mathrm{Ch}_{g,1};\mathbb{Q})
\end{array}$$

$\mathrm{id}_{\bigwedge^2 U_\mathbb{Q}}$

By Theorem 4.3 and this commutative diagram, the kernel

$$\mathrm{Ker}\big((\tau_{g,1}(1))^* : H^2(U;\mathbb{Q}) \to H^2(\mathrm{Ch}_{g,1};\mathbb{Q})\big)$$

contains $[0]_{\mathrm{Sp}} \oplus [2^2]_{\mathrm{Sp}}$, and taking the dual of this, we obtain

$$\mathrm{Im}\big((\tau_{g,1}(1))_* : H_2(\mathrm{Ch}_{g,1};\mathbb{Q}) \to H_2(U;\mathbb{Q})\big) \subset \begin{cases} [1^2]_{\mathrm{Sp}} \oplus [2^2 1^2]_{\mathrm{Sp}} \oplus [1^4]_{\mathrm{Sp}} \oplus [1^6]_{\mathrm{Sp}} & (g \geq 6), \\ [1^2]_{\mathrm{Sp}} \oplus [2^2 1^2]_{\mathrm{Sp}} \oplus [1^4]_{\mathrm{Sp}} & (g = 5), \\ [1^2]_{\mathrm{Sp}} \oplus [2^2 1^2]_{\mathrm{Sp}} & (g = 4), \\ \{0\} & (g = 3). \end{cases}$$

In fact, the summand $[1^2]_{\mathrm{Sp}}$ is not contained in the image

$$\mathrm{Im}\big((\tau_{g,1}(1))_* : H_2(\mathrm{Ch}_{g,1};\mathbb{Q}) \to H_2(U;\mathbb{Q})\big).$$

In this subsection, we show that any other summands except for $[0]_{\mathrm{Sp}}$, $[2^2]_{\mathrm{Sp}}$, and $[1^2]_{\mathrm{Sp}}$ are contained in the image $\mathrm{Im}\big((\tau_{g,1}(1))_* : H_2(\mathrm{Ch}_{g,1};\mathbb{Q}) \to H_2(U;\mathbb{Q})\big)$.

Now, we introduce some $\mathrm{Sp}(2g,\mathbb{Q})$-equivariant homomorphisms to detect specific irreducible component, and *abelian cycles* (see [39; 40]).

Let $V$ be a representation of $\mathrm{Sp}(2g,\mathbb{Q})$.

(1) **The contraction** For $k \geq 2$, $C_k : \bigwedge^k H_\mathbb{Q} \to \bigwedge^{k-2} H_\mathbb{Q}$ is defined by

$$x_1 \wedge \cdots \wedge x_k \mapsto \sum_{i<j} (-1)^{i+j+1} (x_i \cdot x_j) x_1 \wedge \cdots \wedge \widehat{x_i} \wedge \cdots \wedge \widehat{x_j} \wedge \cdots \wedge x_k.$$

Also, the kernel of the contraction $\mathrm{Ker}(C_k)$ corresponds to an irreducible representation denoted by $[1^k]_{\mathrm{Sp}}$.

(2) **The canonical inclusion** $i_V^k : \bigwedge^k V \hookrightarrow \bigotimes^k V$ is defined by

$$v_1 \wedge \cdots \wedge v_k \mapsto \sum_{\sigma \in \mathfrak{S}_k} \mathrm{sign}(\sigma) v_{\sigma(1)} \otimes \cdots \otimes v_{\sigma(k)}.$$

(3) **The multiplication** $\phi_V^{m,n} : (\bigwedge^m V) \otimes (\bigwedge^n V) \to \bigwedge^{m+n} V$ is defined by

$$(v_1 \wedge \cdots \wedge v_m) \otimes (v_{m+1} \wedge \cdots \wedge v_{m+n}) \mapsto v_1 \wedge \cdots \wedge v_m \wedge v_{m+1} \wedge \cdots \wedge v_{m+n}.$$

(4) **The Jacobi identity map** $j_V : \bigwedge^3 V \to V \otimes \bigwedge^2 V$ is defined by

$$v_1 \wedge v_2 \wedge v_3 \mapsto v_1 \otimes (v_2 \wedge v_3) + v_2 \otimes (v_3 \wedge v_1) + v_3 \otimes (v_1 \wedge v_2).$$

Next, we introduce *abelian cycles* which give concrete elements of the image

$$\mathrm{Im}((\tau_{g_1}(1))_* : H_2(\mathrm{Ch}_{g,1}; \mathbb{Z}) \to H_2(U; \mathbb{Z}))$$

of the homomorphism between the second rational homology induced by the first Johnson homomorphism.

**Definition 4.5** [39, Subsection 4.3; 40, page 103, Step 2] Let $G$ be a group and $c : \mathbb{Z}^2 \to G$ be a homomorphism. The image of the fundamental class $1 \in H_2(\mathbb{Z}^2; \mathbb{Z})$ under the induced homomorphism $c_* : H_2(\mathbb{Z}^2; \mathbb{Z}) \to H_2(G; \mathbb{Z})$ is called the *abelian cycle*.

Let $\{e_1, e_2\}$ denote the standard basis of $\mathbb{Z}^2$. Recall that for a finitely generated free abelian group $A$, the second homology group $H_2(A; \mathbb{Z})$ is naturally isomorphic to the second exterior power $\bigwedge^2 A$.

**Proposition 4.6** [37, Lemma 2.1; 39, Lemma 4.5] *Let $A$ be a finitely generated free abelian group and $c : \mathbb{Z}^2 \to A$ be a homomorphism. Then the abelian cycle with respect to $c$ coincides with*

$$c(e_1) \wedge c(e_2) \in \bigwedge^2 A \cong H_2(A; \mathbb{Z}).$$

If we apply this to $\mathbb{Z}^2 \xrightarrow{c} \mathrm{Ch}_{g,1} \xrightarrow{\tau_{g,1}(1)} U$ where $c(e_i) = f_i \in \mathrm{Ch}_{g,1}$ for $i = 1, 2$, we obtain the following.

**Proposition 4.7** *Let $f_1$ and $f_2$ be mutually commutative elements in $\mathrm{Ch}_{g,1}$. Then the element*

$$\tau_{g,1}(1)(f_1) \wedge \tau_{g,1}(1)(f_2) \in \bigwedge^2 U_{\mathbb{Q}} \cong H_2(U; \mathbb{Q})$$

*belongs to the image*

$$\mathrm{Im}((\tau_{g_1}(1))_* : H_2(\mathrm{Ch}_{g,1}; \mathbb{Q}) \to H_2(U; \mathbb{Q})).$$

Additionally, we introduce elements of $\mathrm{Sp}(2g, \mathbb{Z}) \subset \mathrm{Sp}(2g, \mathbb{Q})$ that will appear several times. Let $I$ denote the identity matrix and for distinct $1 \le i, j \le g$, we define the matrix $A_{i,j}$ by the transformation

$$A_{i,j} := \begin{cases} a_i \mapsto a_i + b_i - b_j, \\ a_j \mapsto a_j + b_j - b_i, \\ a_k \mapsto a_k \quad (k \neq i, j), \\ b_k \mapsto b_k. \end{cases}$$

**Proposition 4.8** *For $g \ge 4$, the summand $[2^2 1^2]_{\mathrm{Sp}}$ is contained in the image*

$$\mathrm{Im}\big((\tau_{g,1}(1))_* : H_2(\mathrm{Ch}_{g,1}; \mathbb{Q}) \to H_2(U; \mathbb{Q})\big).$$

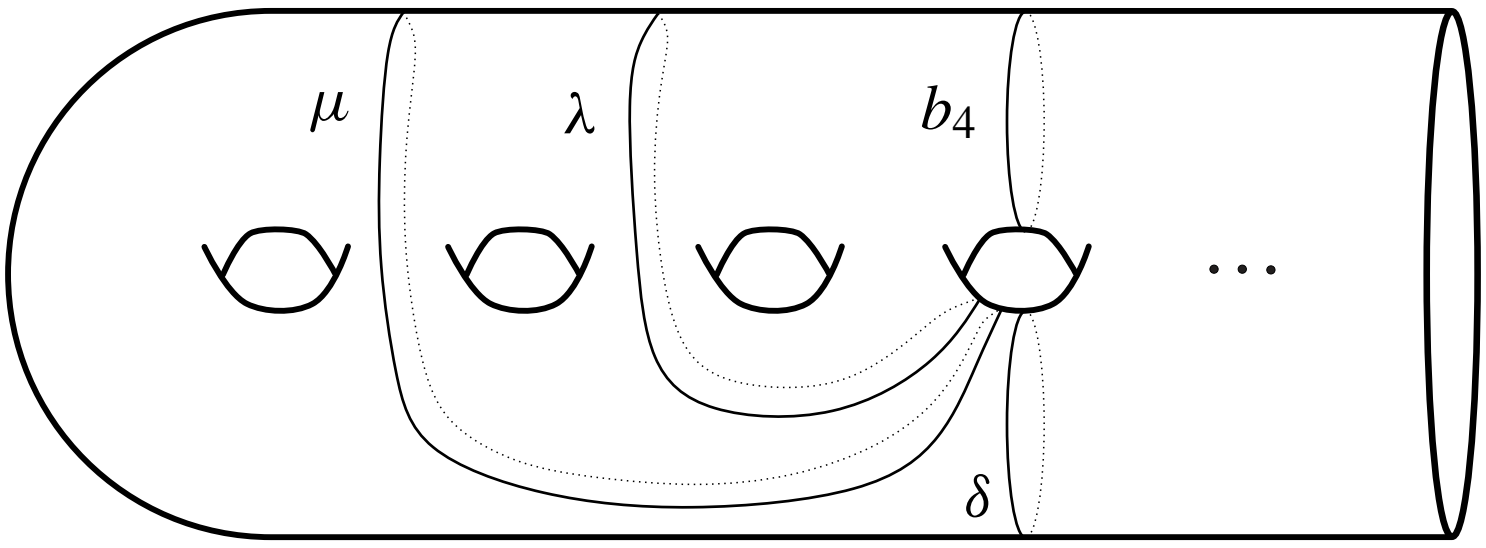

Figure 7: Some simple closed curves on the surface defining an abelian cycle which detects the summand $[2^2 1^2]_{\mathrm{Sp}}$.

**Proof** We take some simple closed curves on the surface as in Figure 7 and we define a homomorphism $\mathbb{Z}^2 \to \mathrm{Ch}_{g,1}$ by

$$\begin{aligned} e_1 &\mapsto \mathrm{BP}(b_4,\delta)\mathrm{BP}(b_4,\mu)^{-1}\mathrm{BP}(b_4,\lambda)^{-1} = T_{b_4}{}^{-1}T_\delta{}^{-1}T_\mu T_\lambda, \\ e_2 &\mapsto \mathrm{BP}(b_4,\mu)\mathrm{BP}(b_4,\lambda)^{-2} = T_{b_4}{}^{-1}T_\mu{}^{-1}T_\lambda{}^2. \end{aligned}$$

We confirm that these two elements are contained in $\mathrm{Ch}_{g,1}$:

$$\begin{aligned} C_3\circ\tau_{g,1}(1)(\mathrm{BP}(b_4,\delta)\,\mathrm{BP}(b_4,\mu)\,\mathrm{BP}(b_4,\lambda)) &= C_3(a_1\wedge b_1\wedge b_4 - a_3\wedge b_3\wedge b_4) = 0, \\ C_3\circ\tau_{g,1}(1)(\mathrm{BP}(b_4,\mu)\mathrm{BP}(b_4,\lambda)^{-2}) &= C_3(a_2\wedge b_2\wedge b_4 - a_3\wedge b_3\wedge b_4) = 0. \end{aligned}$$

Therefore, we obtain

$$\begin{aligned} \zeta_1 &:= (a_1\wedge b_1\wedge b_4 - a_3\wedge b_3\wedge b_4)\wedge(a_2\wedge b_2\wedge b_4 - a_3\wedge b_3\wedge b_4) \\ &= \begin{pmatrix} (a_1\wedge b_1\wedge b_4)\wedge(a_2\wedge b_2\wedge b_4) + (a_2\wedge b_2\wedge b_4)\wedge(a_3\wedge b_3\wedge b_4) \\ +(a_3\wedge b_3\wedge b_4)\wedge(a_1\wedge b_1\wedge b_4) \end{pmatrix} \in \textstyle\bigwedge^2 U_{\mathbb{Q}} \end{aligned}$$

as an element of $\mathrm{Im}\big((\tau_{g,1}(1))_* : H_2(\mathrm{Ch}_{g,1};\mathbb{Q}) \to H_2(U;\mathbb{Q})\big)$.

To prove Proposition 4.8, it is enough to show that $\zeta_1$ is nontrivial on the summand $[2^2 1^2]_{\mathrm{Sp}}$. We detect the nontriviality by using an $\mathrm{Sp}(2g,\mathbb{Q})$-equivariant homomorphism

$$G_1 : \textstyle\bigwedge^2 U_{\mathbb{Q}} \to \big(\bigwedge^4 H_{\mathbb{Q}}\big)\otimes\big(\bigwedge^2 H_{\mathbb{Q}}\big) \supset [2^2 1^2]_{\mathrm{Sp}}$$

as the composition of the maps

$$\begin{aligned} \textstyle\bigwedge^2 U_{\mathbb{Q}} &\hookrightarrow \textstyle\bigwedge^2\big(\bigwedge^3 H_{\mathbb{Q}}\big), \\ i^2_{\bigwedge^3 H_{\mathbb{Q}}} : \textstyle\bigwedge^2\big(\bigwedge^3 H_{\mathbb{Q}}\big) &\to \textstyle\bigotimes^2\big(\bigwedge^3 H_{\mathbb{Q}}\big), \\ \mathrm{id}_{\bigwedge^3 H_{\mathbb{Q}}}\otimes j_{H_{\mathbb{Q}}} : \big(\textstyle\bigwedge^3 H_{\mathbb{Q}}\big)\otimes\big(\bigwedge^3 H_{\mathbb{Q}}\big) &\to \big(\textstyle\bigwedge^3 H_{\mathbb{Q}}\big)\otimes H_{\mathbb{Q}}\otimes\big(\bigwedge^2 H_{\mathbb{Q}}\big), \\ \phi^{3,1}_{H_{\mathbb{Q}}}\otimes\mathrm{id}_{\bigwedge^2 H_{\mathbb{Q}}} : \big(\textstyle\bigwedge^3 H_{\mathbb{Q}}\big)\otimes H_{\mathbb{Q}}\otimes\big(\bigwedge^2 H_{\mathbb{Q}}\big) &\to \big(\textstyle\bigwedge^4 H_{\mathbb{Q}}\big)\otimes\big(\bigwedge^2 H_{\mathbb{Q}}\big). \end{aligned}$$

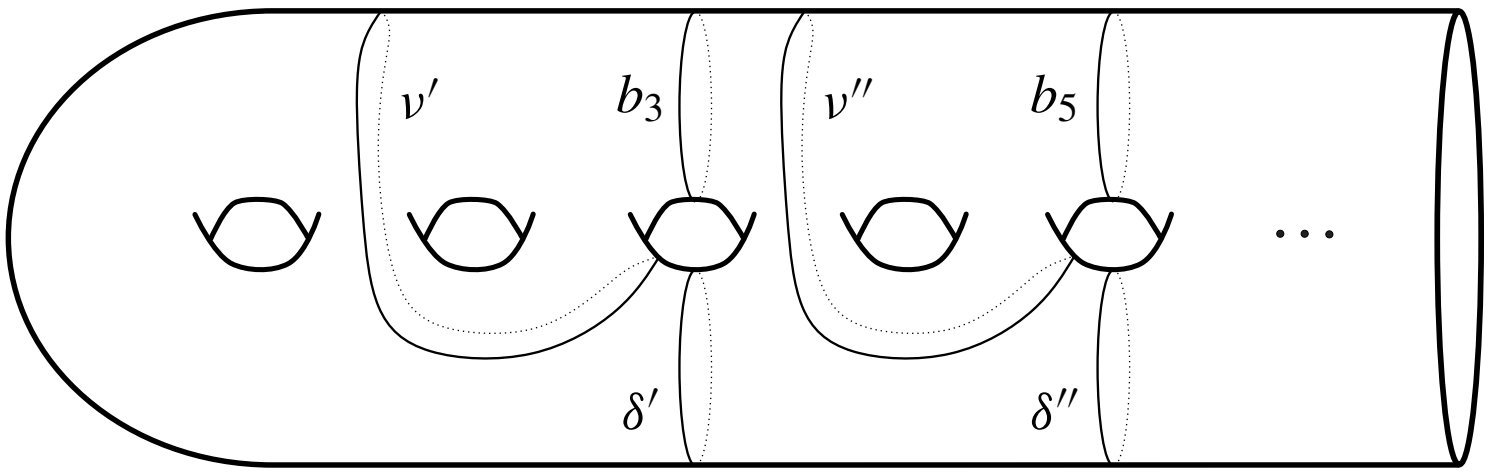

Figure 8: Some simple closed curves on the surface defining an abelian cycle which detects the summand $[1^4]_{\mathrm{Sp}}$.

Using this homomorphism and appropriate elements of $\mathrm{Sp}(2g,\mathbb{Q})$, we compute

$$\begin{aligned}
&\begin{pmatrix}(a_1\wedge b_1\wedge b_4)\wedge(a_2\wedge b_2\wedge b_4)+(a_2\wedge b_2\wedge b_4)\wedge(a_3\wedge b_3\wedge b_4)\\ \qquad +(a_3\wedge b_3\wedge b_4)\wedge(a_1\wedge b_1\wedge b_4)\end{pmatrix}\\
&\quad\xmapsto{I-A_{2,3}}\begin{pmatrix}-2(a_1\wedge b_1\wedge b_4)\wedge(b_2\wedge b_3\wedge b_4)+(a_2\wedge b_2\wedge b_4)\wedge(b_2\wedge b_3\wedge b_4)\\ \qquad +(a_3\wedge b_3\wedge b_4)\wedge(b_2\wedge b_3\wedge b_4)\end{pmatrix}\\
&\quad\xmapsto{I-A_{1,2}} 3(b_1\wedge b_2\wedge b_4)\wedge(b_2\wedge b_3\wedge b_4)\\
&\quad\xmapsto{i^2_{\bigwedge^3 H_{\mathbb{Q}}}} 3(b_1\wedge b_2\wedge b_4)\otimes(b_2\wedge b_3\wedge b_4)-3(b_2\wedge b_3\wedge b_4)\otimes(b_1\wedge b_2\wedge b_4)\\
&\quad\xmapsto{\mathrm{id}_{\bigwedge^3 H_{\mathbb{Q}}}\otimes j_{H_{\mathbb{Q}}}}\begin{pmatrix}3(b_1\wedge b_2\wedge b_4)\otimes(b_2\otimes(b_3\wedge b_4)+b_3\otimes(b_4\wedge b_2)+b_4\otimes(b_2\wedge b_3))\\ -3(b_2\wedge b_3\wedge b_4)\otimes(b_1\otimes(b_2\wedge b_4)+b_2\otimes(b_4\wedge b_1)+b_4\otimes(b_1\wedge b_2))\end{pmatrix}\\
&\quad\xmapsto{\phi^{3,1}_{H_{\mathbb{Q}}}\otimes\mathrm{id}_{\bigwedge^2 H_{\mathbb{Q}}}} -6(b_4\wedge b_2\wedge b_1\wedge b_3)\otimes(b_4\wedge b_2)\\
&\quad\xmapsto[b_4\mapsto b_1,\,b_1\mapsto b_3,\,b_3\mapsto b_4]{a_4\mapsto a_1,\,a_1\mapsto a_3,\,a_3\mapsto a_4,} -6(b_1\wedge b_2\wedge b_3\wedge b_4)\otimes(b_1\wedge b_2)\\
&\quad\xmapsto{b_i\mapsto a_i,\,a_i\mapsto -b_i\ (i=1,2,3,4)} -6(a_1\wedge a_2\wedge a_3\wedge a_4)\otimes(a_1\wedge a_2).
\end{aligned}$$

This vector $-6(a_1\wedge a_2\wedge a_3\wedge a_4)\otimes(a_1\wedge a_2)$ is a highest weight vector of $\left(\bigwedge^4 H_{\mathbb{Q}}\right)\otimes\left(\bigwedge^2 H_{\mathbb{Q}}\right)$. Hence $G_1(\zeta_1)$ is nontrivial on the summand $[2^2 1^2]_{\mathrm{Sp}}$, and $[2^2 1^2]_{\mathrm{Sp}}$ is contained in the image

$$\mathrm{Im}\big((\tau_{g,1}(1))_*: H_2(\mathrm{Ch}_{g,1};\mathbb{Q})\to H_2(U;\mathbb{Q})\big).$$ □

**Proposition 4.9** *For $g\geq 5$, the summand $[1^4]_{\mathrm{Sp}}$ is contained in the image*

$$\mathrm{Im}\big((\tau_{g,1}(1))_*: H_2(\mathrm{Ch}_{g,1};\mathbb{Q})\to H_2(U;\mathbb{Q})\big).$$

**Proof** We take some simple closed curves on the surface as in Figure 8 and we define a homomorphism $\mathbb{Z}^2\to\mathrm{Ch}_{g,1}$ and an abelian cycle as

$$e_1\mapsto \mathrm{BP}(b_3,\delta')\,\mathrm{BP}(b_3,\nu')^{-2}={T_{b_3}}^{-1}{T_{\delta'}}^{-1}{T_{\nu'}}^{2},\quad e_2\mapsto \mathrm{BP}(b_5,\delta'')\,\mathrm{BP}(b_5,\nu'')^{-4}={T_{b_5}}^{-3}{T_{\delta''}}^{-1}{T_{\nu''}}^{4}.$$

Similarly, we obtain

$$\zeta_2 := (a_1 \wedge b_1 \wedge b_3 - a_2 \wedge b_2 \wedge b_3) \wedge (a_1 \wedge b_1 \wedge b_5 + a_2 \wedge b_2 \wedge b_5 + a_3 \wedge b_3 \wedge b_5 - 3 + a_4 \wedge b_4 \wedge b_5)$$

as an element of the image $\operatorname{Im}\big((\tau_{g,1}(1))_* : H_2(\mathrm{Ch}_{g,1};\mathbb{Q}) \to H_2(U;\mathbb{Q})\big)$.

To prove Proposition 4.9, it is enough to show that $\zeta_2$ is nontrivial on the summand $[1^4]_{\mathrm{Sp}}$, and we detect the nontriviality by using an $\mathrm{Sp}(2g,\mathbb{Q})$-equivariant homomorphism $G_2 : \bigwedge^2 U_{\mathbb{Q}} \to \bigwedge^4 H_{\mathbb{Q}} \supset [1^4]_{\mathrm{Sp}}$ as the composition of the maps

$$\begin{aligned} \textstyle\bigwedge^2 U_{\mathbb{Q}} &\hookrightarrow \textstyle\bigwedge^2(\bigwedge^3 H_{\mathbb{Q}}), \\ i^2_{\bigwedge^3 H_{\mathbb{Q}}} : \textstyle\bigwedge^2(\bigwedge^3 H_{\mathbb{Q}}) &\to \textstyle\bigotimes^2(\bigwedge^3 H_{\mathbb{Q}}), \\ \phi_{H^{3,3}_{\mathbb{Q}}} : (\textstyle\bigwedge^3 H_{\mathbb{Q}}) \otimes (\bigwedge^3 H_{\mathbb{Q}}) &\to \textstyle\bigwedge^6 H_{\mathbb{Q}}, \\ C_6 : \textstyle\bigwedge^6 H_{\mathbb{Q}} &\to \textstyle\bigwedge^4 H_{\mathbb{Q}}. \end{aligned}$$

Similarly, we compute

$$\zeta_2 = \begin{pmatrix} (a_1 \wedge b_1 \wedge b_3) \wedge (a_1 \wedge b_1 \wedge b_5) + (a_1 \wedge b_1 \wedge b_3) \wedge (a_2 \wedge b_2 \wedge b_5) \\ + (a_1 \wedge b_1 \wedge b_3) \wedge (a_3 \wedge b_3 \wedge b_5) \\ - 3(a_1 \wedge b_1 \wedge b_3) \wedge (a_4 \wedge b_4 \wedge b_5) - (a_2 \wedge b_2 \wedge b_3) \wedge (a_1 \wedge b_1 \wedge b_5) \\ - (a_2 \wedge b_2 \wedge b_3) \wedge (a_2 \wedge b_2 \wedge b_5) \\ - (a_2 \wedge b_2 \wedge b_3) \wedge (a_3 \wedge b_3 \wedge b_5) + 3(a_2 \wedge b_2 \wedge b_3) \wedge (a_4 \wedge b_4 \wedge b_5) \end{pmatrix}$$

$$\xmapsto{i^2_{\bigwedge^3 H_{\mathbb{Q}}}} \begin{pmatrix} (a_1 \wedge b_1 \wedge b_3) \otimes (a_1 \wedge b_1 \wedge b_5) - (a_1 \wedge b_1 \wedge b_5) \otimes (a_1 \wedge b_1 \wedge b_3) \\ + (a_1 \wedge b_1 \wedge b_3) \otimes (a_2 \wedge b_2 \wedge b_5) - (a_2 \wedge b_2 \wedge b_5) \otimes (a_1 \wedge b_1 \wedge b_3) \\ + (a_1 \wedge b_1 \wedge b_3) \otimes (a_3 \wedge b_3 \wedge b_5) - (a_3 \wedge b_3 \wedge b_5) \otimes (a_1 \wedge b_1 \wedge b_3) \\ - 3(a_1 \wedge b_1 \wedge b_3) \otimes (a_4 \wedge b_4 \wedge b_5) + 3(a_4 \wedge b_4 \wedge b_5) \otimes (a_1 \wedge b_1 \wedge b_3) \\ - (a_2 \wedge b_2 \wedge b_3) \otimes (a_1 \wedge b_1 \wedge b_5) + (a_1 \wedge b_1 \wedge b_5) \otimes (a_2 \wedge b_2 \wedge b_3) \\ - (a_2 \wedge b_2 \wedge b_3) \otimes (a_2 \wedge b_2 \wedge b_5) + (a_2 \wedge b_2 \wedge b_5) \otimes (a_2 \wedge b_2 \wedge b_3) \\ - (a_2 \wedge b_2 \wedge b_3) \otimes (a_3 \wedge b_3 \wedge b_5) + (a_3 \wedge b_3 \wedge b_3) \otimes (a_2 \wedge b_2 \wedge b_3) \\ + 3(a_2 \wedge b_2 \wedge b_3) \otimes (a_4 \wedge b_4 \wedge b_5) - 3(a_4 \wedge b_4 \wedge b_5) \otimes (a_2 \wedge b_2 \wedge b_3) \end{pmatrix}$$

$$\xmapsto{\phi_{H^{3,3}_{\mathbb{Q}}}} 6a_2 \wedge b_2 \wedge a_4 \wedge b_4 \wedge b_3 \wedge b_5 - 6a_1 \wedge b_1 \wedge a_4 \wedge b_4 \wedge b_3 \wedge b_5$$

$$\xmapsto{C_6} 6a_2 \wedge b_2 \wedge b_3 \wedge b_5 - 6a_1 \wedge b_1 \wedge b_3 \wedge b_5 \ (\neq 0)$$

$$\xmapsto{C_4} 6b_3 \wedge b_5 - 6b_3 \wedge b_5 = 0.$$

Since this abelian cycle is nontrivial on the kernel $\operatorname{Ker}(C_4) = [1^4]_{\mathrm{Sp}}$, it follows that $G_2(\zeta_2)$ is nontrivial on the summand $[1^4]_{\mathrm{Sp}}$, and $[1^4]_{\mathrm{Sp}}$ is contained in the image $\operatorname{Im}\big((\tau_{g,1}(1))* : H_2(\mathrm{Ch}_{g,1};\mathbb{Q}) \to H_2(U;\mathbb{Q})\big)$. □

**Proposition 4.10** *For $g \geq 6$, the summand $[1^6]_{\mathrm{Sp}}$ is contained in the image*

$$\mathrm{Im}\big((\tau_{g,1}(1))_* : H_2(\mathrm{Ch}_{g,1};\mathbb{Q}) \to H_2(U;\mathbb{Q})\big).$$

**Proof** For $g \geq 6$, the same abelian cycle as in Proposition 4.9 is also nontrivial on the summand $[1^6]_{\mathrm{Sp}}$, and we check this by using an $\mathrm{Sp}(2g,\mathbb{Q})$-equivariant homomorphism $G_3 : \bigwedge^2 U_{\mathbb{Q}} \to \bigwedge^6 H_{\mathbb{Q}} \supset [1^6]_{\mathrm{Sp}}$ as the composition of the maps

$$\bigwedge^2 U_{\mathbb{Q}} \hookrightarrow \bigwedge^2\big(\bigwedge^3 H_{\mathbb{Q}}\big),$$
$$i^2_{\bigwedge^3 H_{\mathbb{Q}}} : \bigwedge^2\big(\bigwedge^3 H_{\mathbb{Q}}\big) \to \bigotimes^2\big(\bigwedge^3 H_{\mathbb{Q}}\big),$$
$$\phi_{H^{3,3}_{\mathbb{Q}}} : \big(\bigwedge^3 H_{\mathbb{Q}}\big) \otimes \big(\bigwedge^3 H_{\mathbb{Q}}\big) \to \bigwedge^6 H_{\mathbb{Q}}.$$

The result is

$$\begin{aligned}
\zeta_3 := \zeta_2 &= (a_1\wedge b_1\wedge b_3 - a_2\wedge b_2\wedge b_3)\wedge(a_1\wedge b_1\wedge b_5 + a_2\wedge b_2\wedge b_5 + a_3\wedge b_3\wedge b_5 - 3a_4\wedge b_4\wedge b_5)\\
&= \big((a_1\wedge b_1\wedge b_3)\wedge(a_1\wedge b_1\wedge b_5) + (a_1\wedge b_1\wedge b_3)\wedge(a_2\wedge b_2\wedge b_5) + (a_1\wedge b_1\wedge b_3)\wedge(a_3\wedge b_3\wedge b_5)\\
&\quad -3(a_1\wedge b_1\wedge b_3)\wedge(a_4\wedge b_4\wedge b_5) - (a_2\wedge b_2\wedge b_3)\wedge(a_1\wedge b_1\wedge b_5) - (a_2\wedge b_2\wedge b_3)\wedge(a_2\wedge b_2\wedge b_5)\\
&\qquad - (a_2\wedge b_2\wedge b_3)\wedge(a_3\wedge b_3\wedge b_5) + 3(a_2\wedge b_2\wedge b_3)\wedge(a_4\wedge b_4\wedge b_5)\big)\\
&\xmapsto{I-A_{4,6}} 3(a_2\wedge b_2\wedge b_3)\wedge(b_6\wedge b_4\wedge b_5) - 3(a_1\wedge b_1\wedge b_3)\wedge(b_6\wedge b_4\wedge b_5)\\
&\xmapsto{I-A_{1,2}} 6(b_1\wedge b_2\wedge b_3)\wedge(b_4\wedge b_5\wedge b_6)\\
&\xmapsto{i^2_{\bigwedge^3 H_{\mathbb{Q}}}} 6(b_1\wedge b_2\wedge b_3)\otimes(b_4\wedge b_5\wedge b_6) - 6(b_4\wedge b_5\wedge b_6)\wedge(b_1\wedge b_2\wedge b_3)\\
&\xmapsto{\phi_{H^{3,3}_{\mathbb{Q}}}} 12 b_1\wedge b_2\wedge b_3\wedge b_4\wedge b_5\wedge b_6\\
&\xmapsto{b_i\mapsto a_i,\ a_i\mapsto -b_i\ (i=1,2,3,4,5,6)} 12 a_1\wedge a_2\wedge a_3\wedge a_4\wedge a_5\wedge a_6.
\end{aligned}$$

This vector $12a_1\wedge a_2\wedge a_3\wedge a_4\wedge a_5\wedge a_6$ is a highest weight vector of $\bigwedge^6 H_{\mathbb{Q}}$. Hence $G_3(\zeta_3)$ is nontrivial on the summand $[1^6]_{\mathrm{Sp}}$, and $[1^6]_{\mathrm{Sp}}$ is contained in the image $\mathrm{Im}((\tau_{g,1}(1))_* : H_2(\mathrm{Ch}_{g,1};\mathbb{Q}) \to H_2(U;\mathbb{Q}))$. □

Propositions 4.8–4.10 and Theorem 4.3 together imply that the summands $[0]_{\mathrm{Sp}}$ and $[2^2]_{\mathrm{Sp}}$ are contained in the kernel

$$\mathrm{Ker}\big((\tau_{g,1}(1))^* : H^2(U;\mathbb{Q}) \to H^2(\mathrm{Ch}_{g,1};\mathbb{Q})\big),$$

whereas the summands $[2^2 1^2]_{\mathrm{Sp}}$, $[1^4]_{\mathrm{Sp}}$ and $[1^6]_{\mathrm{Sp}}$ are not contained in it. Next, for $g \geq 4$, we prove that the summand $[1^2]_{\mathrm{Sp}}$ is not contained in the image

$$\mathrm{Im}\big((\tau_{g,1}(1))_* : H_2(\mathrm{Ch}_{g,1};\mathbb{Q}) \to H_2(U;\mathbb{Q})\big)$$

and is contained in the kernel

$$\mathrm{Ker}\big((\tau_{g,1}(1))^* : H^2(U;\mathbb{Q}) \to H^2(\mathrm{Ch}_{g,1};\mathbb{Q})\big).$$

**Proposition 4.11** *The diagram*

$$\begin{array}{ccccccc}
H_2(\mathrm{Ch}_{g,1};\mathbb{Q}) & \xrightarrow[(-)^{ab}_*]{} & H_2(\mathrm{Ch}^{ab}_{g,1};\mathbb{Q}) \cong \bigwedge^2 H_1(\mathrm{Ch}_{g,1};\mathbb{Q}) & \xrightarrow{\text{bracket}} & (\Gamma_2(\mathrm{Ch}_{g,1})/\Gamma_3(\mathrm{Ch}_{g,1}))\otimes\mathbb{Q} & \longrightarrow & 0 \\
\downarrow{\scriptstyle (\tau_{g,1}(1))_*} & & \downarrow{\scriptstyle \bigwedge^2(\tau_{g,1}(1))_*} & & \downarrow & & \\
H_2(U;\mathbb{Q}) & \xrightarrow[\cong]{} & \bigwedge^2 H_1(U;\mathbb{Q}) & \xrightarrow{\text{bracket}} & (\mathcal{K}_{g,1}/\mathcal{M}_{g,1}[4])\otimes\mathbb{Q} \cong \mathrm{Im}(\tau_{g,1}(2))\otimes\mathbb{Q} & & \\
 & & \downarrow{\scriptstyle \bigwedge^2 \eta_{\mathbb{Q}}^{-1}} & & \downarrow{\scriptstyle \eta_{\mathbb{Q}}^{-1}} & & \\
 & & \bigwedge^2 \mathcal{T}_1(H_{\mathbb{Q}}) & \xrightarrow{[\bullet,\bullet]_{\mathcal{T}}} & \mathcal{T}_2(H_{\mathbb{Q}}) & & \\
 & & & & \downarrow{\scriptstyle q} & & \\
 & & & \overset{s}{\dashrightarrow} & \bigwedge^2 H_{\mathbb{Q}} & &
\end{array}$$

*is $\mathcal{M}_{g,1}$-equivariant and commutative, where the first row is exact, and the $\mathrm{Sp}(2g,\mathbb{Q})$-equivariant homomorphism $q:\mathcal{T}_2(H_{\mathbb{Q}})\to\bigwedge^2 H_{\mathbb{Q}}$ (see [31]) is defined by*

$q\Bigg($ 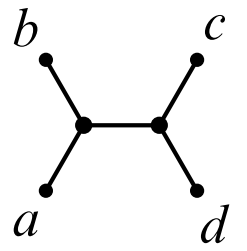 $\Bigg)$

$$:= 4(a\cdot b)(c\wedge d)+4(c\cdot d)(a\wedge b) + 2(d\cdot a)(b\wedge c)+2(b\cdot c)(d\wedge a)+2(a\cdot c)(b\wedge d)+2(d\cdot b)(c\wedge a).$$

**Proof** For the exactness of the first row of the diagram, see [16, page 24, diagram 1.11]. The commutativity in the upper left follows from the naturality of homology with respect to group homomorphisms, while the commutativity in the lower right follows from $\eta_{\mathbb{Q}}^{-1}:\mathfrak{h}_{g,1\,\mathbb{Q}}\to\mathcal{T}(H_{\mathbb{Q}})$ being a homomorphism of Lie algebras. Regarding the commutativity in the upper right, it arises from the brackets induced by commutators within the mapping class group, and the vertical natural homeomorphisms. Recall that

$$\Gamma_2 Ch_{g,1}\subset\Gamma_2\mathcal{I}_{g,1}=\Gamma_2\mathcal{M}_{g,1}[2]\subset\mathcal{M}_{g,1}[3]=\mathcal{K}_{g,1} \quad\text{and}\quad \Gamma_3 Ch_{g,1}\subset\Gamma_3\mathcal{I}_{g,1}=\Gamma_3\mathcal{M}_{g,1}[2]\subset\mathcal{M}_{g,1}[4],$$

due to the fact $\{\mathcal{M}_{g,1}[i+1]/\mathcal{M}_{g,1}[i+2]\}_{i\ge1}=\{\mathrm{Im}(\tau_{g,1}(i))\}_{i\ge1}$ forms a graded Lie algebra under the commutator:

$$\begin{array}{ccc}
\bigwedge^2(\mathrm{Ch}_{g,1}/\Gamma_2(\mathrm{Ch}_{g,1})) & \xrightarrow{\text{bracket}} & \Gamma_2(\mathrm{Ch}_{g,1})/\Gamma_3(\mathrm{Ch}_{g,1}) \\
\downarrow & & \downarrow \\
\bigwedge^2(\mathrm{Ch}_{g,1}/\mathcal{M}_{g,1}[3]) & \xrightarrow{\text{bracket}} & \mathcal{M}_{g,1}[3]/\mathcal{M}_{g,1}[4]
\end{array}$$

□

**Proposition 4.12** *For $g\ge4$, the summand $[1^2]_{\mathrm{Sp}}\subset H_2(U;\mathbb{Q})\cong\bigwedge^2 U_{\mathbb{Q}}$ does not appear in the image*

$$\mathrm{Im}\big((\tau_{g,1}(1))_*:H_2(\mathrm{Ch}_{g,1};\mathbb{Q})\to H_2(U;\mathbb{Q})\big),$$

*and the summand appears in the kernel*

$$\mathrm{Ker}\big((\tau_{g,1}(1))^*:H^2(U;\mathbb{Q})\to H^2(\mathrm{Ch}_{g,1};\mathbb{Q})\big).$$

**Proof** If $[1^2]_{\mathrm{Sp}} \subset H_2(U;\mathbb{Q}) \cong \bigwedge^2 U_{\mathbb{Q}}$ appears in the image of $(\tau_{g,1}(1))_* : H_2(\mathrm{Ch}_{g,1};\mathbb{Q}) \to H_2(U;\mathbb{Q})$, then the $\mathrm{Sp}(2g,\mathbb{Q})$-equivariant homomorphism $s : H_2(U;\mathbb{Q}) \to \bigwedge^2 H_{\mathbb{Q}}$ has to be trivial on $[1^2]_{\mathrm{Sp}}$ because of the commutativity of the diagram and the exactness of the first row. Let

$$\xi_0 = (a_1 \wedge a_3 \wedge b_3 - a_1 \wedge a_4 \wedge b_4) \wedge (a_2 \wedge a_3 \wedge b_3 - a_2 \wedge a_4 \wedge b_4)$$

be an element of $\bigwedge^2 U_{\mathbb{Q}} \cong H_2(U;\mathbb{Q})$. We compute the value of $\xi_0$ under $s$ as

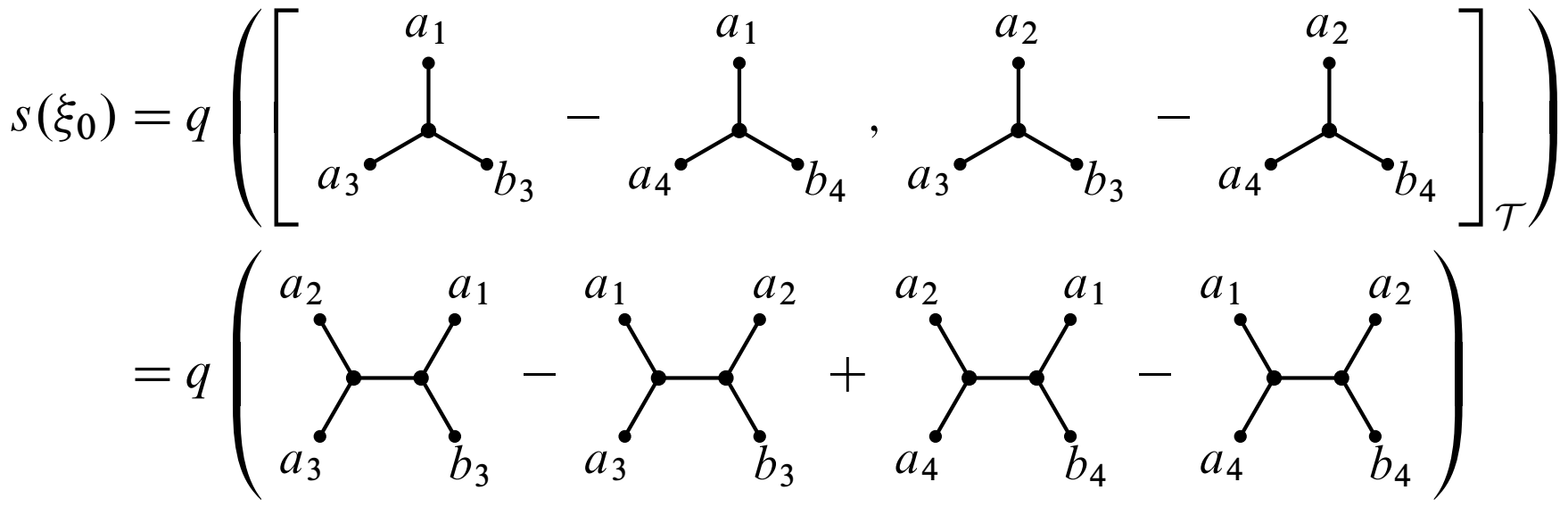

$$= \big(2(b_3 \cdot a_3)(a_2 \wedge a_1) - 2(b_3 \cdot a_3)(a_1 \wedge a_2) + 2(b_4 \cdot a_4)(a_2 \wedge a_1) - 2(b_4 \cdot a_4)(a_1 \wedge a_2)\big)$$
$$= 2a_1 \wedge a_2 + 2a_1 \wedge a_2 + 2a_1 \wedge a_2 + 2a_1 \wedge a_2$$
$$= 8a_1 \wedge a_2.$$

The vector $8a_1 \wedge a_2$ is a highest weight vector of $\bigwedge^2 H_{\mathbb{Q}}$. Hence the $\mathrm{Sp}(2g,\mathbb{Q})$-equivariant homomorphism $s : H_2(U;\mathbb{Q}) \to \bigwedge^2 H_{\mathbb{Q}}$ is nontrivial on the summand $[1^2]_{\mathrm{Sp}}$, which leads to a contradiction. Therefore, the summand $[1^2]_{\mathrm{Sp}} \subset H_2(U;\mathbb{Q})$ never appears in the image of $(\tau_{g,1}(1))_* : H_2(\mathrm{Ch}_{g,1};\mathbb{Q}) \to H_2(U;\mathbb{Q})$ and it does appear in the kernel of $(\tau_{g,1}(1))^* : H^2(U;\mathbb{Q}) \to H^2(\mathrm{Ch}_{g,1};\mathbb{Q})$ for $g \geq 4$. □

From the above considerations, we conclude:

**Theorem 4.13** (Theorem B) *For $g \geq 3$, we have*

$$\mathrm{Im}\big((\tau_{g,1}(1))_* : H_2(\mathrm{Ch}_{g,1};\mathbb{Q}) \to H_2(U;\mathbb{Q})\big) = \begin{cases} [2^2 1^2]_{\mathrm{Sp}} \oplus [1^4]_{\mathrm{Sp}} \oplus [1^6]_{\mathrm{Sp}} & (g \geq 6), \\ [2^2 1^2]_{\mathrm{Sp}} \oplus [1^4]_{\mathrm{Sp}} & (g = 5), \\ [2^2 1^2]_{\mathrm{Sp}} & (g = 4), \\ \{0\} & (g = 3), \end{cases}$$

*and*

$$\mathrm{Ker}\big((\tau_{g,1}(1))^* : H^2(U;\mathbb{Q}) \to H^2(\mathrm{Ch}_{g,1};\mathbb{Q})\big) = \begin{cases} [0]_{\mathrm{Sp}} \oplus [2^2]_{\mathrm{Sp}} \oplus [1^2]_{\mathrm{Sp}} & (g \geq 4), \\ [0]_{\mathrm{Sp}} \oplus [2^2]_{\mathrm{Sp}} & (g = 3) \end{cases}$$

*as* $\mathrm{Sp}(2g,\mathbb{Q})$*-modules, and the same holds for the* $\mathrm{Ch}_{g,*}$ *case.*

Moreover, for a 2-cocycle $v$ representing an element of the kernel

$$\mathrm{Ker}\big((\tau_{g,1}(1))^* : H^2(U;\mathbb{Q}) \to H^2(\mathrm{Ch}_{g,1};\mathbb{Q})\big) \cong (\tau_{g,1}(2) \otimes \mathbb{Q})^*,$$

we explicitly construct a coboundary of $(\tau_{g,1}(1))^*(v)$.

**Proposition 4.14** *Let* $[v] \in \mathrm{Ker}(H^2(U, \mathbb{Q}) \to H^2(\mathrm{Ch}_{g,1}; \mathbb{Q})) \cong (\mathrm{Im}(\tau_{g,1}(2)) \otimes \mathbb{Q})^* \cong \mathcal{T}_2(H_{\mathbb{Q}})^*$, *where* $v$ *is a representative* 2*-cocycle. Then the* 1*-cochain* $v \circ 2r_2^\theta$ *cobounds* $(\tau_{g,1}(1))^*(v)$.

**Proof** We identify $\mathrm{Ker}(H^2(U; \mathbb{Q}) \to H^2(\mathrm{Ch}_{g,1}; \mathbb{Q}))$ with the image of the map

$$[\bullet, \bullet]^* : \mathcal{T}_2(H_{\mathbb{Q}})^* \cong (\mathrm{Im}(\tau_{g,1}(2)) \otimes \mathbb{Q})^* \hookrightarrow H^2(U; \mathbb{Q})$$

dual of the Lie bracket. For $[v] \in \mathrm{Ker}(H^2(U; \mathbb{Q}) \to H^2(\mathrm{Ch}_{g,1}; \mathbb{Q}))$, the 2-cocycle

$$(\tau_{g,1}(1))^*(v) = v \circ [\tau_{g,1}, \tau_{g,1}] = v \circ [r_1^\theta, r_1^\theta]_{\mathcal{T}}$$

is cobounded by $v \circ r_2^\theta$. In fact, from the second-degree part of the BCH series, we compute

$$\begin{aligned}\delta(v \circ 2r_2^\theta)(f, g) &= v \circ 2(r_2^\theta(fg) - r_2^\theta(f) - r_2^\theta(g)) \\ &= v \circ 2\big(\big(\tfrac{1}{2}[r_1^\theta(f), r_1^\theta(g)]_{\mathcal{T}} + r_2^\theta(f) + r_2^\theta(g)\big) - r_2^\theta(f) - r_2^\theta(g)\big) \\ &= v \circ 2\big(\tfrac{1}{2}[r_1^\theta(f), r_1^\theta(g)]_{\mathcal{T}}\big) \\ &= v \circ [r_1^\theta(f), r_1^\theta(g)]_{\mathcal{T}}.\end{aligned}$$

□

# 5 The Casson–Morita homomorphism $d : \mathrm{Ch}_{g,1} \to \mathbb{Z}$ for the Chillingworth subgroup

Morita [25] introduced a certain *map* $d : \mathcal{M}_{g,1} \to \mathbb{Z}$ related to the Casson invariant.

**Definition 5.1** A Heegaard embedding $\Sigma_{g,1} \to \Sigma_g \to S^3$ is a map such that cutting along it decomposes into two genus-$g$ handlebodies $V_g^+$ and $V_g^-$.

**Definition 5.2** For an element $\varphi \in \mathcal{I}_{g,1}$, the integral homology 3-sphere $M_\varphi$ is defined as the 3-manifold obtained by regluing $V_g^-$ to $V_g^+$ along $\varphi$.

The Casson invariant $\lambda$ : {integral homology 3-spheres} $\to \mathbb{Z}$ is one of the fundamental invariants of integral homology 3-spheres. In particular, we can consider the Casson invariant for $M_\varphi$. He found that the Casson invariant can be interpreted as a secondary invariant associated with the characteristic classes of the surface and studied this mapping in detail. Morita showed in [28, Proposition 2.1] that the map $\varphi \mapsto \lambda(M_\varphi) =: \lambda^*(\varphi)$ (the map $\lambda^*$ depends on the choice of a Heegaard embedding and there is no canonical choice) is a homomorphism on the Johnson kernel $\mathcal{K}_{g,1}$. He defined a homomorphism $d = d|_{\mathcal{K}_{g,1}} : \mathcal{K}_{g,1} \to \mathbb{Z}$ related to $\lambda^*$. He called this homomorphism the core of the Casson invariant, and we call it the *Casson–Morita homomorphism*.

To define the Casson–Morita homomorphism, we introduce some 2-cocycles of the full mapping class group $\mathcal{M}_{g,1}$.

**Definition 5.3** Let $\tau : \mathcal{M}_{g,1} \times \mathcal{M}_{g,1} \to \mathbb{Z}$ be the *Meyer cocycle* characterized by the signature of the 4-manifold defined by the surface $\Sigma_g$ bundle over a pair of pants $\Sigma_{0,3}$ with corresponding monodromies (see [24]). Next, let $k : \mathcal{M}_{g,1} \to H^{(*)}$ be a crossed homomorphism representing a generator of $H^1(\mathcal{M}_{g,1}; H^{(*)}) \cong \mathbb{Z}$, for example the Chillingworth homomorphism $k = e_X$. We define the 2-cocycle $c : \mathcal{M}_{g,1} \times \mathcal{M}_{g,1} \to \mathbb{Z}$ by $c(\varphi, \psi) := k(\varphi^{-1}) \cdot k(\psi)$ called the *intersection cocycle*.

These 2-cocycles are related by $[-3\tau] = e_1 = [c] \in H^2(\mathcal{M}_{g,1};\mathbb{Z})$, where $e_1$ is the first *Mumford–Morita–Miller class* (see [24; 25; 27]). Therefore, there exists a map $d : \mathcal{M}_{g,1} \to \mathbb{Z}$ such that the coboundary $\delta d$ coincides with $c + 3\tau$ as 2-cocycles. Moreover, for $g \geq 3$, $H^1(\mathcal{M}_{g,1};\mathbb{Z}) = 0$ holds (Mumford [35], Birman [3] and Powell [36] showed this for the closed case. For the general case, see a Korkmaz's survey [21]). Hence, such a map $d : \mathcal{M}_{g,1} \to \mathbb{Z}$ is uniquely determined. Therefore, we will always assume $g \geq 3$ from now on.

We have the following by definition.

**Proposition 5.4** *For $f, g \in \mathcal{M}_{g,1}$, we have*

$$d(fg) = d(f) + d(g) - k(f^{-1}) \cdot k(g) - 3\tau(f, g).$$

**Definition 5.5** We define the homomorphism $d = d|_{\mathrm{Ch}_{g,1}} : \mathrm{Ch}_{g,1} \to \mathbb{Z}$ given above as the Casson–Morita homomorphism on $\mathrm{Ch}_{g,1}$.

By this equality, $d = d|_{\mathrm{Ch}_{g,1}} : \mathrm{Ch}_{g,1} \to \mathbb{Z}$ is a homomorphism on the Chillingworth subgroup because the Meyer cocycle $\tau$ is realized as a normalized 2-cocycle on $\mathrm{Sp}(2g,\mathbb{Z})$, hence vanishes on the Torelli group $\mathcal{I}_{g,1}$ (see [24]), and the crossed homomorphism $k$ is trivial on the Chillingworth subgroup $\mathrm{Ch}_{g,1}$.

**Remark 5.6** The Casson–Morita *map* $d : \mathcal{M}_{g,1} \to \mathbb{Z}$ depends on the choice of a crossed homomorphism $k : \mathcal{M}_{g,1} \to H^{(*)}$.

**Proposition 5.7** *The Casson–Morita homomorphism $d : \mathrm{Ch}_{g,1} \to \mathbb{Z}$ does not depend on the choice of a crossed homomorphism $k : \mathcal{M}_{g,1} \to H^*$.*

**Proof** Let $k_0, k_1 : \mathcal{M}_{g,1} \to H^*$ be two crossed homomorphisms representing a generator of the group $H^1(\mathcal{M}_{g,1}; H^*) \cong \mathbb{Z}$. For $i = 0, 1$, we denote the intersection cocycles determined by $k_i$ as $c_i$ and the Casson–Morita homomorphisms as $d_i$. First, when $k_1 = -k_0$, since $c_0 = c_1$, we have $d_0 = d_1$ on $\mathcal{M}_{g,1}$. Therefore, it suffices to consider the case where $k_1$ is cohomologous to $k_0$ in $H^1(\mathcal{M}_{g,1}; H^*)$. In this case, we can write $k_1(f) - k_0(f) = (f^{-1})^*h - h$ for some element $h \in H^*$. Then, $d_1(f) - d_0(f) = (k_1(f^{-1}) - k_0(f)) \cdot h$, and the right side is always 0 on the Chillingworth subgroup $\mathrm{Ch}_{g,1}$. Indeed, the calculation proceeds by direct computation as follows:

$$\begin{aligned}
\delta(d_1 - d_0)(f,g) &= (c_1 - c_0)(f,g) \\
&= k_1(f^{-1}) \cdot k_1(g) - k_0(f^{-1}) \cdot k_0(g) \\
&= (k_0(f^{-1}) + f^*h - h) \cdot (k_0(g) + (g^{-1})^*h - h) - k_0(f^{-1}) \cdot k_0(g) \\
&= k_0(f^{-1}) \cdot ((g^{-1})^*h - h) + (f^*h - h) \cdot k_0(g) + (f^*h - h) \cdot ((g^{-1})^*h - h) \\
&= (g^*k_0(f^{-1}) - k_0(f^{-1})) \cdot h - ((f^{-1})^*k_0(g) - k_0(g)) \cdot h \\
&\qquad\qquad + (g^*f^*h - f^*h - g^*h) \cdot h \\
&= (k_0((fg)^{-1}) - k_0(f^{-1}) - k_0(g^{-1})) \cdot h - (k_0(fg) - k_0(f) - k_0(g)) \cdot h \\
&\qquad\qquad + ((fg)^*h - f^*h - g^*h) \cdot h \\
&= \delta\big(\bullet \mapsto (k_0(\bullet^{-1}) + (\bullet)^*h - h - k_0(\bullet)) \cdot h\big)(f,g) \\
&= \delta\big(\bullet \mapsto (k_1(\bullet^{-1}) - k_0(\bullet)) \cdot h\big)(f,g).
\end{aligned}$$

Since $H^1(\mathcal{M}_{g,1};\mathbb{Z})=0$, the coboundary is unique, leading to

$$d_1(f)-d_0(f)=(k_1(f^{-1})-k_0(f))\cdot h,$$

and hence $d_1$ and $d_0$ coincide on $\mathrm{Ch}_{g,1}$. □

Morita [25] gave some properties and formulas of the Casson–Morita map.

**Proposition 5.8** (Morita [25, page 320, Theorem 5.3]) (1) *The Casson–Morita homomorphism $d:\mathcal{M}_{g,1}\to\mathbb{Z}$ defined by a crossed homomorphism $k:\mathcal{M}_{g,1}\to H$ is stable with respect to the genus of the surface if $k$ is stable with respect to the genus. Specifically, for the homomorphism $\mathcal{M}_{g,1}\to\mathcal{M}_{g+1,1}$ induced by gluing a genus-1 surface with 2 boundary components, denoted by $\Sigma_{1,2}$, to the boundary of $\Sigma_{g,1}$ to obtain $\Sigma_{g+1,1}$, the diagram*

$$\begin{array}{ccc} \mathcal{M}_{g,1} & & \\ \downarrow & \overset{d}{\searrow} & \\ \mathcal{M}_{g+1,1} & \xrightarrow{d} & \mathbb{Z} \end{array}$$

*commutes if the diagram*

$$\begin{array}{ccc} \mathcal{M}_{g,1} & \xrightarrow{k} & H_1(\Sigma_{g,1};\mathbb{Z}) \\ \downarrow & & \downarrow \\ \mathcal{M}_{g+1,1} & \xrightarrow{k} & H_1(\Sigma_{g+1,1};\mathbb{Z}) \end{array}$$

*commutes.*

(2) *Let $T_\gamma$ be a genus-$h$ BSCC map, meaning that $\gamma$ bounds a genus-$h$ subsurface of the surface. Then the value under $d$ is $4h(h-1)$. In particular, the Casson–Morita homomorphism $d:\mathcal{K}_{g,1}\to\mathbb{Z}$ on the Johnson kernel is $\mathcal{M}_{g,1}$-invariant, and its image coincides with $8\mathbb{Z}$.*

**Proposition 5.9** (part of Theorem C) *The Casson–Morita map $d:\mathrm{Ch}_{g,1}\to\mathbb{Z}$ on the Chillingworth subgroup is also an $\mathcal{M}_{g,1}$-invariant homomorphism.*

**Proof** Consider $h\in\mathrm{Ch}_{g,1}$ and $f\in\mathcal{M}_{g,1}$. Then, we have

$$\begin{aligned} d(fhf^{-1}) &= d(fh)+d(f^{-1})-k((fh)^{-1})\cdot k(f^{-1})-3\tau(fh,f^{-1}) \\ &= (d(f)+d(h)-k(f^{-1})\cdot k(h)-3\tau(f,h)) \\ &\qquad -d(f)-(k(h^{-1})+h^*k(f^{-1}))\cdot k(f^{-1})-3\tau(f,f^{-1}) \\ &= d(f)+d(h)-d(f)-3\tau(f,\mathrm{id}_{\Sigma_{g,1}})-k(f^{-1})\cdot k(f^{-1})-3\tau(f,f^{-1}) \\ &= d(h), \end{aligned}$$

where we used following properties of the Meyer cocycle $\tau:\mathcal{M}_{g,1}\times\mathcal{M}_{g,1}\to\mathbb{Z}$ (see [24]).

(1) *The Meyer cocycle factors through $\mathrm{Sp}(2g,\mathbb{Z})$, which means that for any $h_1$, $h_2\in\mathcal{I}_{g,1}$, we have $\tau(fh_1,gh_2)=\tau(f,g)$.*

(2) $\tau(f,f^{-1})=0$.

(3) $\tau(f,\mathrm{id}_{\Sigma_{g,1}})=0$. □

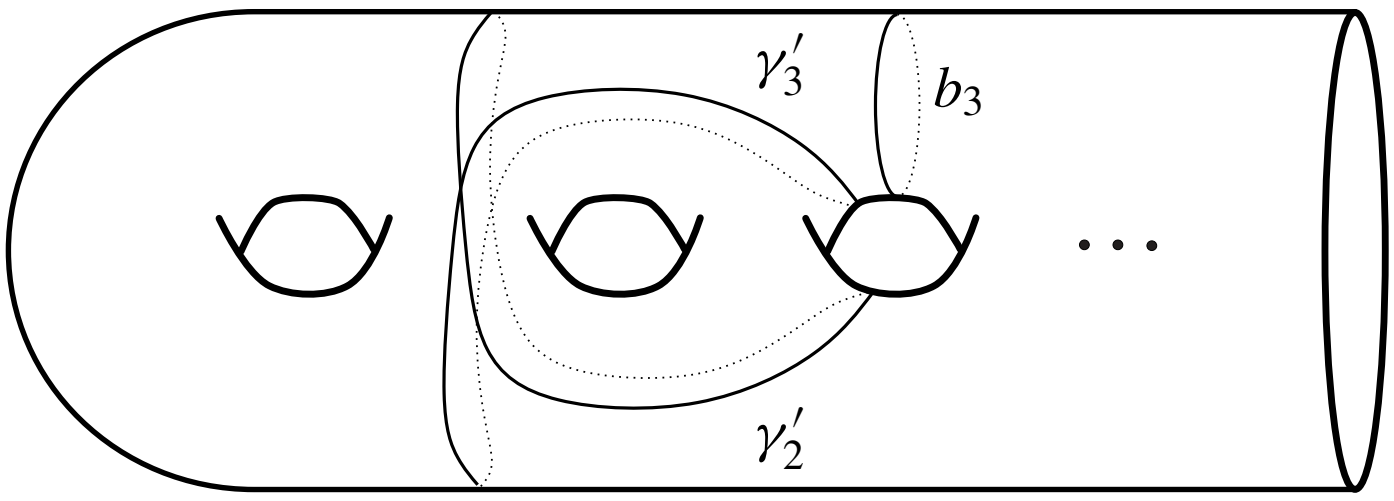

Figure 9: Simple closed curves $\gamma_2'$, $\gamma_3'$ defining a homological genus-zero bounding pair map $B_0 := T_{\gamma_2'} T_{\gamma_3'}{}^{-1}$ and $b_3$.

**Proposition 5.10** (part of Theorem C) *The image of the Casson–Morita homomorphism $d : \mathrm{Ch}_{g,1} \to \mathbb{Z}$ on the Chillingworth subgroup is also $8\mathbb{Z}$.*

**Proof** Let $k$ denote the Chillingworth class $t$. By

$$d(\text{genus-}h \text{ BSCC map}) = 4h(h-1),$$

the image of the Casson–Morita homomorphism restricted to the Johnson kernel $\mathcal{K}_{g,1}$ is $8\mathbb{Z}$. Let us consider the element $B_0 := \mathrm{BP}(\gamma_2', \gamma_3') := T_{\gamma_2'} T_{\gamma_3'}{}^{-1}$ in Figure 9. We have $d(B_0) = 0$. Indeed, we consider $k(\mathrm{BP}(\gamma_2', b_3)) - k(\mathrm{BP}(\gamma_3', b_3)) = 2[\gamma_2'] - 2[\gamma_3'] = 0$ and using $k(fg) = k(f) + f_* k(g)$, we have $k(T_{\gamma_2'}) = k(T_{\gamma_3'})$, and note that $\rho(T_{\gamma_2'}) = \rho(T_{\gamma_3'}) \in \mathrm{Sp}(2g, \mathbb{Z})$. Moreover, considering the braid relations

$$T_{a_3} T_{\gamma_i'} T_{a_3} = T_{\gamma_i'} T_{a_3} T_{\gamma_i'}$$

and using $k(T_{\gamma_2'}) = k(T_{\gamma_3'})$ and $\rho(T_{\gamma_2'}) = \rho(T_{\gamma_3'})$, we have $d(T_{\gamma_2'}) = d(T_{\gamma_3'}) = -d(T_{\gamma_3'}^{-1})$ and consequently, $d(B_0) = 0$. Since $\mathrm{Ch}_{g,1}$ is normally generated by the Johnson kernel $\mathcal{K}_{g,1}$ and $B_0$, the image of the map $d : \mathrm{Ch}_{g,1} \to \mathbb{Z}$ coincides with $8\mathbb{Z}$. □

We also determine the kernel of the Casson–Morita homomorphism for the Chillingworth subgroup. Before discussing it, we present Faes's result on the Johnson kernel, which provides the motivation for our study.

**Theorem 5.11** (Faes [11, Remark 2.15]) *For $g \geq 2$, the kernel of the Casson–Morita homomorphism restricted to the Johnson kernel is given by*

$$\mathrm{Ker}(d|_{\mathcal{K}_{g,1}} : \mathcal{K}_{g,1} \to \mathbb{Z}) = \langle T_{\gamma_1'} \rangle [\mathcal{K}_{g,1}, \mathcal{M}_{g,1}],$$

*where $T_{\gamma_1'}$ is the Dehn twist along $\gamma_1'$ also called a genus-one BSCC map as shown in Figure 10, $\langle T_{\gamma_1'} \rangle$ is the subgroup generated by $T_{\gamma_1'}$, and $[\mathcal{K}_{g,1}, \mathcal{M}_{g,1}]$ is the commutator subgroup of the Johnson kernel and the full mapping class group.*

**Proof sketch** This theorem is essentially based on the result $H^1(\mathcal{K}_{g,1}; \mathbb{Z})^{\mathcal{M}_{g,1}} \cong \mathbb{Z} \oplus \mathbb{Z}$ by Morita [28], where the superscript means the $\mathcal{M}_{g,1}$-invariants. Morita introduced an $\mathcal{M}_{g,1}$-invariant homomorphism $d' : \mathcal{K}_{g,1} \to \mathbb{Z}$, distinct from $d$, and showed that $H^1(\mathcal{K}_{g,1};\mathbb{Q})^{\mathcal{M}_{g,1}}$ is generated by $d$ and $d'$. Faes [11] took linear combinations $\frac{d}{8}$ and $\frac{4d'-5d}{12}$, and proved that these two elements give an isomorphism

$$\overline{\left(\frac{d}{8}, \frac{4d'-5d}{12}\right)} : \mathcal{K}_{g,1}/[\mathcal{K}_{g,1}, \mathcal{M}_{g,1}] \cong H_1(\mathcal{K},\mathbb{Z})_{\mathcal{M}_{g,1}} \xrightarrow{\cong} \mathbb{Z} \oplus \mathbb{Z},$$

and we have $\left(\frac{d}{8}, \frac{4d'-5d}{12}\right)(T_{\gamma_1'}) = (0,1)$. In particular, the intersection of these kernels

$$\operatorname{Ker}\left(\left(\frac{d}{8}, \frac{4d'-5d}{12}\right) : \mathcal{K}_{g,1} \to \mathbb{Z} \oplus \mathbb{Z}\right) = \operatorname{Ker}(d|_{\mathcal{K}_{g,1}}) \cap \operatorname{Ker}(d'|_{\mathcal{K}_{g,1}})$$

coincides with the commutator subgroup $[\mathcal{K}_{g,1}, \mathcal{M}_{g,1}]$ of the Johnson kernel and the full mapping class group. For any elements $\varphi \in \operatorname{Ker}(d|_{\mathcal{K}_{g,1}})$, the element

$$T_{\gamma_1'}^{\left(-\frac{4d'(\varphi)-5d(\varphi)}{12}\right)}\varphi = T_{\gamma_1'}^{\left(-\frac{d'(\varphi)}{3}\right)}\varphi$$

is contained in

$$\operatorname{Ker}(d|_{\mathcal{K}_{g,1}}) \cap \operatorname{Ker}(d'|_{\mathcal{K}_{g,1}}) = [\mathcal{K}_{g,1}, \mathcal{M}_{g,1}].$$

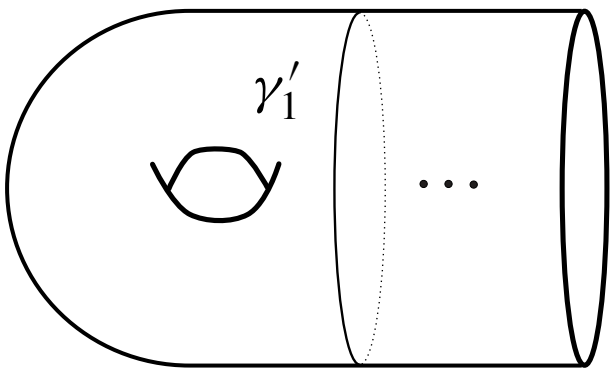

Figure 10: A simple closed curve $\gamma_1'$ on the surface.

Therefore, we have

$$\operatorname{Ker}(d|_{\mathcal{K}_{g,1}}) = \langle T_{\gamma_1'} \rangle [\mathcal{K}_{g,1}, \mathcal{M}_{g,1}]. \qquad \square$$

**Theorem 5.12** (Theorem C) *For $g \geq 4$, the kernel of the Casson–Morita homomorphism on the Chillingworth subgroup is given by*

$$\operatorname{Ker}(d : \mathrm{Ch}_{g,1} \to \mathbb{Z}) = \langle\langle B_0 \rangle\rangle \langle T_{\gamma_1'} \rangle [\mathcal{K}_{g,1}, \mathcal{M}_{g,1}],$$

*where the element $B_0 := T_{\gamma_2'} T_{\gamma_3'}^{-1}$ is a homological genus-zero BP map as shown in Figure 11, and $\langle\langle B_0 \rangle\rangle$ is the normal subgroup of $\mathcal{M}_{g,1}$ generated by $B_0$, and where $T_{\gamma_1'}$ is the Dehn twist along the simple closed curve $\gamma_1'$ in Figure 10.*

**Proof** The element $B_0$ satisfies $d(B_0) = 0$ and $U$ is generated by

$$\tau_{g,1}(1)(B_0) = a_1 \wedge b_1 \wedge b_3 - a_2 \wedge b_2 \wedge b_3$$

as an $\mathrm{Sp}(2g,\mathbb{Z})$-module (see the proof of Proposition 3.14). Therefore, for any $\varphi \in \operatorname{Ker}(d : \mathrm{Ch}_{g,1} \to \mathbb{Z})$, there exists an $\psi \in \langle\langle B_0 \rangle\rangle$ such that $\tau_{g,1}(1)(\psi) = \tau_{g,1}(1)(\varphi)$, which implies that $\psi^{-1}\varphi \in \operatorname{Ker}(d|_{\mathcal{K}_{g,1}})$. Hence, we have

$$\operatorname{Ker}(d : \mathrm{Ch}_{g,1} \to \mathbb{Z}) = \langle\langle B_0 \rangle\rangle \langle T_{\gamma_1'} \rangle [\mathcal{K}_{g,1}, \mathcal{M}_{g,1}]. \qquad \square$$

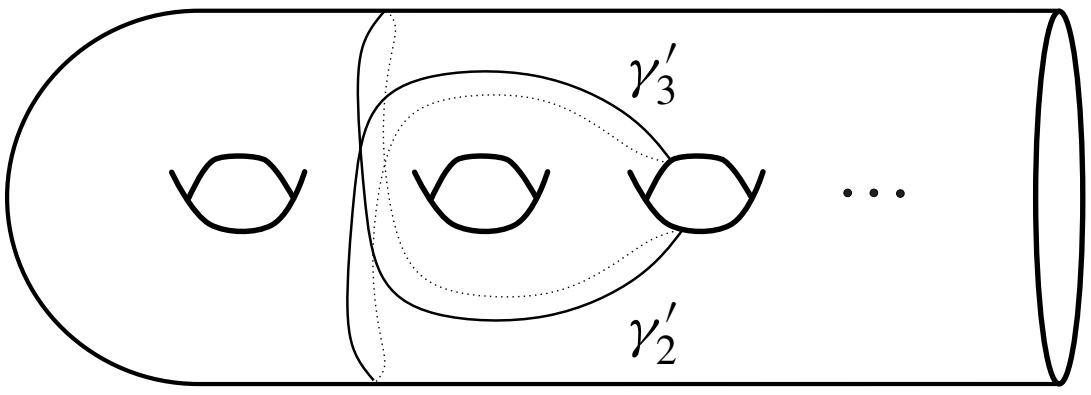

Figure 11: Some simple closed curves on the surface defining a homological genus-zero BP map.

## 6 Proof of Theorem A

For the Torelli group, the rational abelianization is obtained from the first Johnson homomorphism as a mapping class group module. More precisely, Johnson [19] showed that the abelianization of the Torelli group is isomorphic to the direct sum of the target space of the Johnson homomorphism and additional 2-torsion parts, arising from the Birman–Craggs homomorphism. The latter is closely related to spin structures and the Rokhlin invariant (see [4]). For the Chillingworth subgroup, the first Johnson homomorphism provides one abelian quotient. In addition, the Casson–Morita homomorphism is a homomorphism on the Chillingworth subgroup that is nontrivial on the kernel of the first Johnson homomorphism. This allows us to combine both homomorphisms to obtain a better lower bound for the rational abelianization $H_1(\mathrm{Ch}_{g,1};\mathbb{Q}) \cong (\mathrm{Ch}_{g,1})^{ab} \otimes \mathbb{Q}$ of the Chillingworth subgroup:

$$d \oplus \tau_{g,1}(1) : \mathrm{Ch}_{g,1} \to (8\mathbb{Z} \oplus U) \otimes \mathbb{Q}.$$

To determine the rational abelianization of the Chillingworth subgroup, we consider the inflation-restriction exact sequence of the rational homology for the short exact sequence

$$1 \to \mathcal{K}_{g,1} \to \mathrm{Ch}_{g,1} \to U \to 0$$

induced by the first Johnson homomorphism for the Chillingworth subgroup as follows:

$$H_2(\mathrm{Ch}_{g,1};\mathbb{Q}) \to H_2(U;\mathbb{Q}) \cong \textstyle\bigwedge^2 U_{\mathbb{Q}} \to H_1(\mathcal{K}_{g,1};\mathbb{Q})_U \to H_1(\mathrm{Ch}_{g,1};\mathbb{Q}) \to H_1(U;\mathbb{Q}) \cong U_{\mathbb{Q}} \to 0.$$

This exact sequence is equivariant under the natural action of the mapping class group. Having already determined the $\mathcal{M}_{g,1}$-module structures of the image

$$\mathrm{Im}((\tau_{g,1}(1))_* : H_2(\mathrm{Ch}_{g,1};\mathbb{Q}) \to H_2(U;\mathbb{Q}))$$

and $U_{\mathbb{Q}}$, we only have to determine the $\mathcal{M}_{g,1}$-module structure of $H_1(\mathcal{K}_{g,1};\mathbb{Q})_U$, where the subscript $U$ means the $U$-coinvariant of the first rational homology group $H_1(\mathcal{K}_{g,1};\mathbb{Q})$ of the Johnson kernel. To study the structure of $H_1(\mathcal{K}_{g,1};\mathbb{Q})_U$, we use the rational abelianization of the Johnson kernel $\mathcal{K}_{g,1}$ by Faes and Massuyeau [12]. The rational abelianization of the Johnson kernel was originally computed for the case of the closed surface $\mathcal{K}_g$ by Dimca–Hain–Papadima [9] and Morita–Sakasai–Suzuki [34].

Here, we introduce a certain homomorphism $\mathrm{Tr}_3$, which is needed to describe the rational abelianization of the Johnson kernel.

**Definition 6.1** (Morita [29, Section 6]) The $\mathrm{Sp}(2g,\mathbb{Q})$-equivariant homomorphism $\mathrm{Tr}_3 : \mathcal{T}_3(H_{\mathbb{Q}}) \to \mathrm{S}^3 H_{\mathbb{Q}}$ is called Morita's trace map and is defined by

$\mathrm{Tr}_3$ 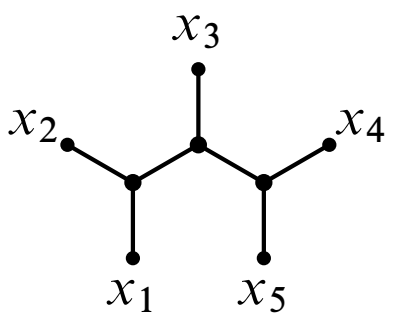

$$:= 2(x_5\cdot x_1)x_2x_3x_4 + 2(x_1\cdot x_4)x_5x_3x_2 + 2(x_4\cdot x_2)x_1x_3x_5 + 2(x_2\cdot x_5)x_2x_3x_1,$$

where $\mathrm{S}^3 H_{\mathbb{Q}}$ is the third symmetric power of $H_{\mathbb{Q}}$.

**Theorem 6.2** (Faes–Massuyeau [12, Theorem 3.2]) *For $g \geq 6$, the rational abelianization of the Johnson kernel $H_1(\mathcal{K}_{g,1};\mathbb{Q})$ as an $\mathcal{M}_{g,1}$-module is given by the Casson–Morita homomorphism $d$ and the truncations of the infinitesimal Dehn–Nielsen representation $(r_2^\theta, r_3^\theta)$ as*

$$d \oplus (r_2^\theta, r_3^\theta) : \mathcal{K}_{g,1} \to \mathbb{Q} \oplus (\mathcal{T}_2(H_{\mathbb{Q}}) \times \mathrm{Ker}(\mathrm{Tr}_3)) \subset \mathbb{Q} \oplus (\mathcal{T}_2(H_{\mathbb{Q}}) \times \mathcal{T}_3(H_{\mathbb{Q}})).$$

*The homomorphism $d \oplus (r_2^\theta, r_3^\theta)$ induces the isomorphism $H_1(\mathcal{K}_{g,1};\mathbb{Q}) \cong \mathbb{Q} \oplus (\mathcal{T}_2(H_{\mathbb{Q}}) \times \mathrm{Ker}(\mathrm{Tr}_3))$ as $\mathcal{M}_{g,1}$-modules.*

**Remark 6.3** The action of the mapping class group $\mathcal{M}_{g,1}$ on $\mathbb{Q} \oplus (\mathcal{T}_2(H_{\mathbb{Q}}) \times \mathrm{Ker}(\mathrm{Tr}_3))$ does not factor through the integral symplectic group $\mathrm{Sp}(2g,\mathbb{Z})$ because of the influence of the bracket.

Next, to study the action of $U \cong \mathrm{Ch}_{g,1}/\mathcal{K}_{g,1}$ on the rational abelianization of the Johnson kernel $H_1(\mathcal{K}_{g,1};\mathbb{Q})$, we summarize the behavior by conjugation.

**Lemma 6.4** *For $f \in \mathcal{I}_{g,1}$ and $h \in \mathcal{K}_{g,1}$, we have*

$$d \oplus (r_2^\theta, r_3^\theta)(fhf^{-1}) = d \oplus (r_2^\theta, r_3^\theta)(h) + \bigl(0, (0, [r_1^\theta(f), r_2^\theta(h)]_{\mathcal{T}})\bigr).$$

**Proof** Since the Casson–Morita homomorphism $d : \mathcal{K}_{g,1} \to \mathbb{Z}$ is $\mathcal{M}_{g,1}$-invariant, it suffices to consider the conjugation action on $(r_2^\theta, r_3^\theta)$. For $f \in \mathcal{I}_{g,1}$ and $h \in \mathcal{K}_{g,1}$, we compute the part of $r^\theta(fhf^{-1})$ up to the third degree, that is, modulo the terms of degree four and higher, $\widehat{\mathcal{T}_{\geq 4}} = \widehat{\bigoplus_{i\geq 4}\mathcal{T}_i(H_{\mathbb{Q}})} \subset \widehat{\mathcal{T}(H_{\mathbb{Q}})}$, as

$$\begin{aligned} r^\theta(fhf^{-1}) &= r^\theta(f) \star r^\theta(hf^{-1}) \\ &\equiv \begin{pmatrix} r^\theta(f) + r^\theta(hf^{-1}) + \frac{1}{2}[r^\theta(f), r^\theta(hf^{-1})]_{\widehat{\mathcal{T}}} \\ + \frac{1}{12}\begin{pmatrix} [r^\theta(f), [r^\theta(f), r^\theta(hf^{-1})]_{\widehat{\mathcal{T}}}]_{\widehat{\mathcal{T}}} \\ \qquad - [r^\theta(hf^{-1}), [r^\theta(f), r^\theta(hf^{-1})]_{\widehat{\mathcal{T}}}]_{\widehat{\mathcal{T}}} \end{pmatrix} \end{pmatrix} \pmod{\widehat{\mathcal{T}_{\geq 4}}}. \end{aligned}$$

Using $r^\theta = r_1^\theta + r_2^\theta + \cdots$, we expand the brackets up to terms of degree 3, we have

$$\begin{aligned} &r^\theta(fhf^{-1}) \\ &\equiv \begin{pmatrix} r^\theta(f) + r^\theta(hf^{-1}) \\ + \frac{1}{2}\Bigl([r_1^\theta(f), r_1^\theta(hf^{-1})]_{\mathcal{T}} + [r_1^\theta(f), r_2^\theta(hf^{-1})]_{\mathcal{T}} + [r_2^\theta(f), r_1^\theta(hf^{-1})]_{\mathcal{T}}\Bigr) \\ + \frac{1}{12}\Bigl([r_1^\theta(f), [r_1^\theta(f), r_1^\theta(hf^{-1})]_{\mathcal{T}}]_{\mathcal{T}} - [r_1^\theta(hf^{-1}), [r_1^\theta(f), r_1^\theta(hf^{-1})]_{\mathcal{T}}]_{\mathcal{T}}\Bigr) \end{pmatrix} \pmod{\widehat{\mathcal{T}_{\geq 4}}}. \end{aligned}$$

By separating the BCH product by each degree, we note that

$$r_1^\theta(hf^{-1}) = r_1^\theta(h) + r_1^\theta(f^{-1}) = r_1^\theta(h) - r_1^\theta(f)$$

and

$$r_2^\theta(hf^{-1}) = r_2^\theta(h) + r_2^\theta(f^{-1}) + \tfrac{1}{2}[r_1^\theta(h), r_1^\theta(f^{-1})] = r_2^\theta(h) - r_2^\theta(f) - \tfrac{1}{2}[r_1^\theta(h), r_1^\theta(f)]_{\mathcal{T}}$$

are obtained. Using this,

$$r^\theta(fhf^{-1}) \equiv \begin{pmatrix} r^\theta(f) + r^\theta(hf^{-1}) + \frac{1}{2}([r_1^\theta(f), r_1^\theta(h) - r_1^\theta(f)]_{\mathcal{T}} \\ \qquad + [r_1^\theta(f), r_2^\theta(h) - r_2^\theta(f) - \frac{1}{2}[r_1^\theta(h), r_1^\theta(f)]_{\mathcal{T}}]_{\mathcal{T}} \\ \qquad + [r_2^\theta(f), r_1^\theta(h) - r_1^\theta(f)]_{\mathcal{T}}) \\ + \frac{1}{12}\begin{pmatrix} [r_1^\theta(f), [r_1^\theta(f), r_1^\theta(h) - r_1^\theta(f)]_{\mathcal{T}}]_{\mathcal{T}} \\ \quad - [r_1^\theta(h) - r_1^\theta(f), [r_1^\theta(f), r_1^\theta(h) - r_1^\theta(f)]_{\mathcal{T}}]_{\mathcal{T}} \end{pmatrix} \end{pmatrix} \pmod{\widehat{\mathcal{T}_{\geq 4}}}.$$

Since $h$ is an element of the Johnson kernel $\mathcal{K}_{g,1}$ and $r_1^\theta$ is the first Johnson homomorphism, we have $r_1^\theta(h) = 0$, which implies

$$\begin{aligned} r^\theta(fhf^{-1}) &\equiv \begin{pmatrix} r^\theta(f) + r^\theta(hf^{-1}) + \frac{1}{2}([r_1^\theta(f), 0 - r_1^\theta(f)]_{\mathcal{T}} \\ \quad + [r_1^\theta(f), r_2^\theta(h) - r_2^\theta(f) - \frac{1}{2}[0, r_1^\theta(f)]_{\mathcal{T}}]_{\mathcal{T}} + [r_2^\theta(f), 0 - r_1^\theta(f)]_{\mathcal{T}}) \\ + \frac{1}{12}\left([r_1^\theta(f), [r_1^\theta(f), 0 - r_1^\theta(f)]_{\mathcal{T}}]_{\mathcal{T}} - [0 - r_1^\theta(f), [r_1^\theta(f), 0 - r_1^\theta(f)]_{\mathcal{T}}]_{\mathcal{T}}\right) \end{pmatrix} \pmod{\widehat{\mathcal{T}_{\geq 4}}} \\ &= r^\theta(f) + r^\theta(hf^{-1}) + \tfrac{1}{2}\left([r_1^\theta(f), r_2^\theta(h) - r_2^\theta(f)]_{\mathcal{T}} - [r_2^\theta(f), r_1^\theta(f)]_{\mathcal{T}}\right) \pmod{\widehat{\mathcal{T}_{\geq 4}}} \\ &= r^\theta(f) + r^\theta(hf^{-1}) + \tfrac{1}{2}[r_1^\theta(f), r_2^\theta(h)]_{\mathcal{T}} \pmod{\widehat{\mathcal{T}_{\geq 4}}}. \end{aligned}$$

Similarly, by expanding the BCH product and collecting nontrivial terms up to the third degree, we have

$$\begin{aligned} r^\theta(fhf^{-1}) &\equiv \begin{pmatrix} r^\theta(f) + r^\theta(h) - r^\theta(f) + \frac{1}{2}[r_1^\theta(f), r_2^\theta(h)]_{\mathcal{T}} \\ \qquad + \frac{1}{2}[r^\theta(h), -r^\theta(f)]_{\widehat{\mathcal{T}}} \\ + \frac{1}{12}\left([r^\theta(h), [r^\theta(h), -r^\theta(f)]_{\widehat{\mathcal{T}}}]_{\widehat{\mathcal{T}}} - [-r^\theta(f), [r^\theta(h), -r^\theta(f)]_{\widehat{\mathcal{T}}}]_{\widehat{\mathcal{T}}}\right) \end{pmatrix} \pmod{\widehat{\mathcal{T}_{\geq 4}}} \\ &\equiv r^\theta(f) + r^\theta(h) - r^\theta(f) + \tfrac{1}{2}[r_1^\theta(f), r_2^\theta(h)]_{\mathcal{T}} - \tfrac{1}{2}[r_2^\theta(h), r_1^\theta(f)]_{\mathcal{T}} \pmod{\widehat{\mathcal{T}_{\geq 4}}} \\ &= r^\theta(h) + [r_1^\theta(f), r_2^\theta(h)]_{\mathcal{T}} \pmod{\widehat{\mathcal{T}_{\geq 4}}}. \end{aligned}$$

Therefore, for $f \in \mathcal{I}_{g,1}$ and $h \in \mathcal{K}_{g,1}$, we have

$$(r_2^\theta, r_3^\theta)(fhf^{-1}) = (r_2^\theta, r_3^\theta)(h) + \left(0, [r_1^\theta(f), r_2^\theta(h)]_{\mathcal{T}}\right).$$ □

**Proposition 6.5** *For* $g \geq 3$, *the* $\mathrm{Sp}(2g, \mathbb{Q})$*-equivariant bracket map*

$$[\bullet, \bullet]_{\mathcal{T}} : (r_1^\theta(\mathrm{Ch}_{g,1}) \otimes \mathbb{Q}) \otimes (r_2^\theta(\mathcal{K}_{g,1}) \otimes \mathbb{Q}) \to \mathrm{Ker}(\mathrm{Tr}_3) \xrightarrow[\cong]{\eta^{\mathbb{Q}}} \mathrm{Im}(\tau_{g,1}(3)) \otimes \mathbb{Q}$$

*is surjective.*

**Proof** The target space

$$\mathrm{Ker}(\mathrm{Tr}_3) \cong \mathrm{Im}(\tau_{g,1}(3)) \otimes \mathbb{Q}$$

(regarding the inclusion $\eta(\mathrm{Ker}(\mathrm{Tr}_3)) \supset \mathrm{Im}(\tau_{g,1}(3))$, see Morita [29, Theorem 6.1]; for the coincidence, see [33]) is isomorphic to $[31^2]_{\mathrm{Sp}} \oplus [21]_{\mathrm{Sp}}$ as representations of the rational symplectic group $\mathrm{Sp}(2g, \mathbb{Q})$, which is shown by Asada and Nakamura [2, Theorem A-(iii)]. Therefore, it is sufficient to show that the bracket map is nontrivial on both the summands $[31^2]_{\mathrm{Sp}}$ and $[21]_{\mathrm{Sp}}$. Let us consider two elements $\xi_1, \xi_2 \in [r_1^\theta(\mathrm{Ch}_{g,1}) \otimes \mathbb{Q}, r_2^\theta(\mathcal{K}_{g,1}) \otimes \mathbb{Q}]_{\mathcal{T}} \subset \mathrm{Ker}(\mathrm{Tr}_3) \subset r_3^\theta(\mathrm{Ch}_{g,1}) \otimes \mathbb{Q} \subset \mathcal{T}_3(H_{\mathbb{Q}})$ defined as

$\xi_1 :=$ 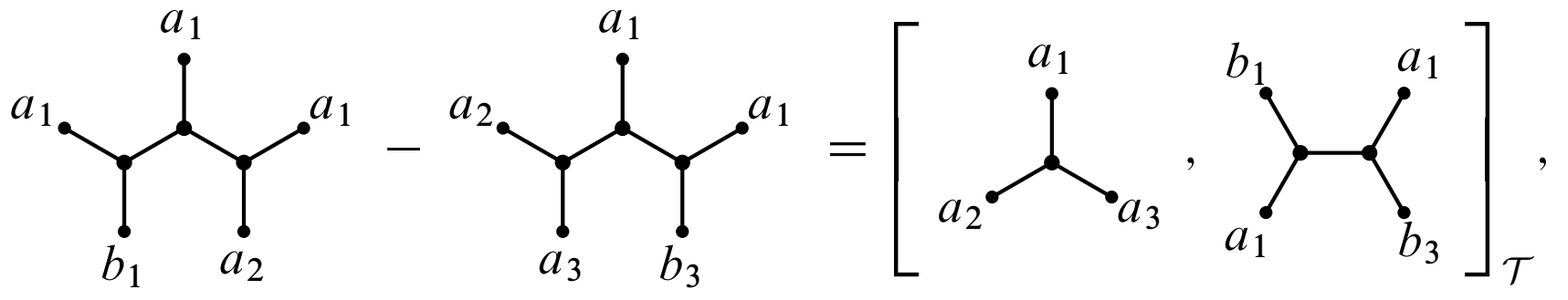

$\xi_2 := -2$ 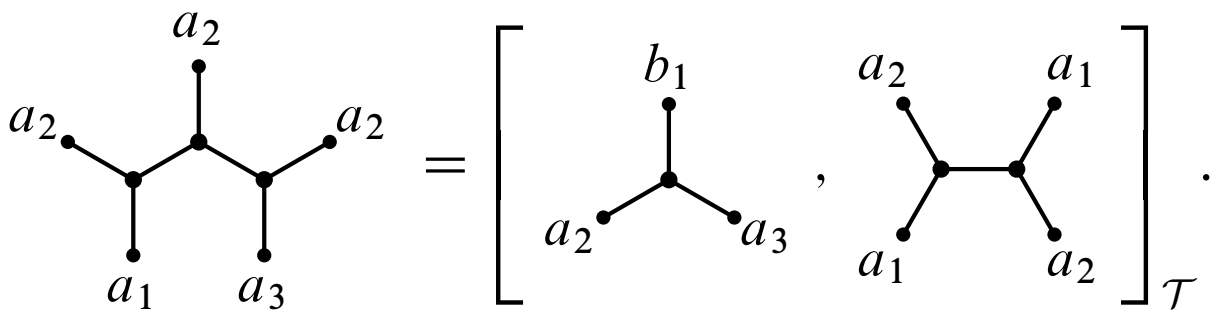

We also define an $\mathrm{Sp}(2g, \mathbb{Q})$-equivariant homomorphism to detect these summands as

$$\mathrm{Ker}(\mathrm{Tr}_3) \hookrightarrow \mathcal{T}_3(H_{\mathbb{Q}}) \hookrightarrow H_{\mathbb{Q}} \otimes \mathcal{L}_{g,1}(4)_{\mathbb{Q}} \hookrightarrow \bigotimes^5 H_{\mathbb{Q}}$$

and

$$\bigotimes^5 H_{\mathbb{Q}} \xrightarrow{x_1 \otimes x_2 \otimes x_3 \otimes x_4 \otimes x_5 \mapsto (x_1 \cdot x_2)(x_3 \wedge x_4) \otimes x_5} \left(\bigwedge^2 H_{\mathbb{Q}}\right) \otimes H_{\mathbb{Q}}.$$

The value of $\xi_1$ under the above homomorphism is $9(a_1 \wedge a_2) \otimes a_1$ which hits the highest weight vector of the summand $[21]_{\mathrm{Sp}}$, and the value of $\xi_2$ under the above homomorphism is 0, but $\xi_2$ itself is nontrivial. Hence, $\xi_2$ purely lies in the summand $[31^2]_{\mathrm{Sp}}$. Therefore, the bracket map is surjective over $\mathbb{Q}$. (We calculated with a Mathematica program based on the method described in [38, p. 22, 7. Appendix].) □

We get the following directly from Lemma 6.4, Proposition 6.5 and the fact, shown by Morita [25, Proposition 1.2] and Hain [15], that $\mathrm{Im}(\tau_{g,1}(2)) \otimes \mathbb{Q} = \mathfrak{h}_{g,1\,\mathbb{Q}}(2) \cong \mathcal{T}_2(H_{\mathbb{Q}}) \cong [2^2]_{\mathrm{Sp}} \oplus [1^2]_{\mathrm{Sp}} \oplus [0]_{\mathrm{Sp}}$.

**Proposition 6.6** *The $U$-coinvariant $H_1(\mathcal{K}_{g,1}; \mathbb{Q})_U$ of the first rational homology of the Johnson kernel are isomorphic to $\mathbb{Q} \oplus \mathcal{T}_2(H_{\mathbb{Q}})$ via the homomorphism $(d, r_2^\theta = \tau_{g,1}(2))$ as $\mathcal{M}_{g,1}$-modules. In particular, the action of the mapping class group on $H_1(\mathcal{K}_{g,1}; \mathbb{Q})_U$ factors through the rational symplectic group $\mathrm{Sp}(2g, \mathbb{Q})$, and decomposes as $[0]_{\mathrm{Sp}} \oplus ([2^2]_{\mathrm{Sp}} \oplus [1^2]_{\mathrm{Sp}} \oplus [0]_{\mathrm{Sp}})$.*

We now proceed to Theorem A.

**Theorem 6.7** (part of Theorem A) *For $g \geq 6$, the first rational (co)homology of the Chillingworth subgroup $\mathrm{Ch}_{g,1}$ for the genus-$g$ surface with one boundary is induced by the Casson–Morita homomorphism and* *the first Johnson homomorphism for the Chillingworth subgroup* $d \oplus \tau_{g,1}(1) : \mathrm{Ch}_{g,1} \to 8\mathbb{Z} \oplus U$, *and satisfies*

$$(\mathrm{Ch}_{g,1})^{ab} \otimes \mathbb{Q} \cong H_1(\mathrm{Ch}_{g,1};\mathbb{Q}) \cong [1^3]_{\mathrm{Sp}} \oplus [0]_{\mathrm{Sp}},$$
$$((\mathrm{Ch}_{g,1})^{ab} \otimes \mathbb{Q})^* \cong H^1(\mathrm{Ch}_{g,1};\mathbb{Q}) \cong [1^3]_{\mathrm{Sp}} \oplus [0]_{\mathrm{Sp}}$$

*as* $\mathcal{M}_{g,1}$*-modules.*

**Proof** Now, we handle the inflation-restriction exact sequence of the rational homology for the short exact sequence

$$1 \to \mathcal{K}_{g,1} \to \mathrm{Ch}_{g,1} \to U \to 0$$

to determine the first rational homology group $H_1(\mathrm{Ch}_{g,1};\mathbb{Q})$ of the Chillingworth subgroup. For $g \geq 6$, we have determined the image

$$\mathrm{Im}\big(H_2(\mathrm{Ch}_{g,1};\mathbb{Q}) \to \textstyle\bigwedge^2 U_{\mathbb{Q}}\big) \cong [1^6]_{\mathrm{Sp}} \oplus [2^2 1^2]_{\mathrm{Sp}} \oplus [1^4]_{\mathrm{Sp}}$$

of the homomorphism between the second rational homology induced by the first Johnson homomorphism for the Chillingworth subgroup, the $U$-coinvariant

$$H_1(\mathcal{K}_{g,1};\mathbb{Q})_U \cong [0]_{\mathrm{Sp}} \oplus ([2^2]_{\mathrm{Sp}} \oplus [1^2]_{\mathrm{Sp}} \oplus [0]_{\mathrm{Sp}})$$

of the first rational homology of the Johnson kernel and $U_{\mathbb{Q}} \cong [1^3]_{\mathrm{Sp}}$. By adding the information obtained from the above to the long exact sequence, we obtain

$$\begin{array}{ccccc}
H_2(\mathrm{Ch}_{g,1};\mathbb{Q}) & \xrightarrow{(\tau_{g,1(1)})^*} & \bigwedge^2 U_{\mathbb{Q}} & \longrightarrow & H_1(\mathcal{K}_{g,1};\mathbb{Q})_U \\
 & & ([1^6]_{\mathrm{Sp}}\oplus[2^2 1^2]_{\mathrm{Sp}}\oplus[1^4]_{\mathrm{Sp}}) & & ([2^2]_{\mathrm{Sp}}\oplus[1^2]_{\mathrm{Sp}}\oplus[0]_{\mathrm{Sp}}) \\
 & & \oplus([2^2]_{\mathrm{Sp}}\oplus[1^2]_{\mathrm{Sp}}\oplus[0]_{\mathrm{Sp}}) & & \oplus[0]_{\mathrm{Sp}} \\
\to H_1(\mathrm{Ch}_{g,1};\mathbb{Q}) & \longrightarrow & U_{\mathbb{Q}} & \longrightarrow & 0. \\
\text{“}[0]_{\mathrm{Sp}} & & [1^3]_{\mathrm{Sp}} & & \\
\oplus[1^3]_{\mathrm{Sp}}\text{”} & & & &
\end{array}$$

The preceding argument alone does not determine whether $H_1(\mathrm{Ch}_{g,1};\mathbb{Q})$ decomposes into a direct sum of two summands, $[0]_{\mathrm{Sp}}$ and $[1^3]_{\mathrm{Sp}}$, as $\mathcal{M}_{g,1}$-modules. However, $\dim_{\mathbb{Q}} H_1(\mathrm{Ch}_{g,1};\mathbb{Q}) = \dim_{\mathbb{Q}}([0]_{\mathrm{Sp}} \oplus [1^3]_{\mathrm{Sp}})$. Combining this with the lower bound of the rational abelianization of the Chillingworth subgroup already obtained, $d \oplus \tau_{g,1}(1) : \mathrm{Ch}_{g,1} \to \mathbb{Q} \oplus U_{\mathbb{Q}} \cong [0]_{\mathrm{Sp}} \oplus [1^3]$ gives a rational abelianization of the Chillingworth subgroup. Therefore, this long exact sequence splits at the $H_1(\mathrm{Ch}_{g,1};\mathbb{Q})$ as $\mathcal{M}_{g,1}$-modules. □

**Corollary 6.8** *For* $g \geq 6$, *the rank of the abelianization of the Chillingworth subgroup* $\mathrm{Ch}_{g,1}$ *for the surface* $\Sigma_{g,1}$ *is* $\frac{1}{3}(2g-1)(2g^2-2g-3)$.

Next, we consider the case of the Chillingworth subgroup $\mathrm{Ch}_{g,*}$ with a fixed base point.

**Theorem 6.9** (part of Theorem A) *For $g \geq 6$, the first rational (co)homology of the Chillingworth subgroup $\mathrm{Ch}_{g,*}$ for the genus-$g$ surface with a base point is induced by the first Johnson homomorphism $\tau_{g,*}(1) : \mathrm{Ch}_{g,*} \to U$, and satisfies*

$$(\mathrm{Ch}_{g,*})^{ab} \otimes \mathbb{Q} \cong H_1(\mathrm{Ch}_{g,*}; \mathbb{Q}) \cong [1^3]_{\mathrm{Sp}}, \quad ((\mathrm{Ch}_{g,*})^{ab} \otimes \mathbb{Q})^* \cong H^1(\mathrm{Ch}_{g,*}; \mathbb{Q}) \cong [1^3]_{\mathrm{Sp}}$$

*as $\mathcal{M}_{g,*}$-modules.*

**Proof** We consider the central extension

$$0 \to \mathbb{Z} \to \mathrm{Ch}_{g,1} \to \mathrm{Ch}_{g,*} \to 1$$

induced by the natural homomorphism $\mathcal{M}_{g,1} \to \mathcal{M}_{g,*}$ and the inflation-restriction exact sequence for it

$$\cdots \to H_1(\mathbb{Z}; \mathbb{Q}) \to H_1(\mathrm{Ch}_{g,1}; \mathbb{Q}) \to H_1(\mathrm{Ch}_{g,*}; \mathbb{Q}) \to 0.$$

The Casson–Morita homomorphism induces the map $\mathbb{Z} \hookrightarrow \mathrm{Ch}_{g,1} \xrightarrow{d} \mathbb{Z}$, $1 \mapsto 4g(g-1)$ which is nontrivial by the formula by Morita. Therefore, the homomorphism $H_1(\mathbb{Z}; \mathbb{Q}) \to H_1(\mathrm{Ch}_{g,1}; \mathbb{Q})$ is nontrivial and the image

$$\mathrm{Im}(H_1(\mathbb{Z}; \mathbb{Q}) \to H_1(\mathrm{Ch}_{g,1}; \mathbb{Q}))$$

coincides with the summand $[0]_{\mathrm{Sp}}$. We have the exact sequence

$$\cdots \to \underset{[0]_{\mathrm{Sp}}}{H_1(\mathbb{Z}; \mathbb{Q})} \to \underset{[0]_{\mathrm{Sp}} \oplus [1^3]_{\mathrm{Sp}}}{H_1(\mathrm{Ch}_{g,1}; \mathbb{Q})} \to \underset{[1^3]_{\mathrm{Sp}}}{H_1(\mathrm{Ch}_{g,*}; \mathbb{Q})} \to 0.$$

In particular, the rational abelianization of the Chillingworth subgroup $\mathrm{Ch}_{g,*}$ for the surface with a base point $\mathrm{Ch}_{g,1}$ is induced by the first Johnson homomorphism for the Chillingworth subgroup alone. □

**Corollary 6.10** *For $g \geq 6$, the rank of the abelianization of the Chillingworth subgroup $\mathrm{Ch}_{g,*}$ for the surface $\Sigma_{g,*}$ is $\frac{2}{3}g(2g+1)(g-2)$.*

Finally, we consider the case of the Chillingworth subgroup $\mathrm{Ch}_g$ without a fixed base point and boundary components.

**Theorem 6.11** (part of Theorem A) *For $g \geq 6$, the first rational (co)homology of the Chillingworth subgroup $\mathrm{Ch}_g$ for the genus-$g$ closed surface $\mathrm{Ch}_g$ is induced by the first Johnson homomorphism $\tau_g(1) : \mathrm{Ch}_g \to \overline{U}$, and satisfies*

$$(\mathrm{Ch}_g)^{ab} \otimes \mathbb{Q} \cong H_1(\mathrm{Ch}_g; \mathbb{Q}) \cong [1^3]_{\mathrm{Sp}}, \quad ((\mathrm{Ch}_g)^{ab} \otimes \mathbb{Q})^* \cong H^1(\mathrm{Ch}_g; \mathbb{Q}) \cong [1^3]_{\mathrm{Sp}}$$

*as $\mathcal{M}_g$-modules.*

**Proof** Let us consider the short exact sequence

$$1 \to [\pi_1(\Sigma_g), \pi_1(\Sigma_g)] \to \mathrm{Ch}_{g,*} \to \mathrm{Ch}_g \to 1$$

induced by the natural homomorphism $\mathcal{M}_{g,*} \to \mathcal{M}_g$ and the long exact sequence for it

$$\cdots \to H_1(\mathrm{Ch}_{g,*};\mathbb{Q}) \to H_1(\mathrm{Ch}_g;\mathbb{Q}) \to 0.$$

Since the rational abelianization of the Chillingworth subgroup $H_1(\mathrm{Ch}_{g,*};\mathbb{Q})$ is irreducible as an $\mathcal{M}_g$-module and there exists the first Johnson homomorphism $\tau_g(1) : \mathrm{Ch}_g \to \overline{U} \subset \bigwedge^3 H/H$, the natural homomorphism $H_1(\mathrm{Ch}_{g,*};\mathbb{Q}) \to H_1(\mathrm{Ch}_g;\mathbb{Q})$ is an isomorphism. □

**Corollary 6.12** *For $g \geq 6$, the rank of the abelianization of the Chillingworth subgroup $\mathrm{Ch}_g$ for the surface $\Sigma_g$ is also $\frac{2}{3}g(2g+1)(g-2)$.*

In the last of this section, we also mention the Euler class of the central extension

$$0 \to \mathbb{Z} \to \mathrm{Ch}_{g,1} \to \mathrm{Ch}_{g,*} \to 1$$

induced by the natural homomorphism $\mathcal{M}_{g,1} \to \mathcal{M}_{g,*}$.

**Theorem 6.13** (Theorem D) *The Euler class $e \in H^2(\mathrm{Ch}_{g,*};\mathbb{Z})$ of the central extension*

$$0 \to \mathbb{Z} \to \mathrm{Ch}_{g,1} \to \mathrm{Ch}_{g,*} \to 1$$

*is a $\frac{g(g-1)}{2}$-torsion element.*

**Proof** We consider the inflation-restriction exact sequence of the integral cohomology

$$\cdots \to H^1(\mathrm{Ch}_{g,1};\mathbb{Z}) \to H^1(\mathbb{Z};\mathbb{Z}) \to H^2(\mathrm{Ch}_{g,*};\mathbb{Z}) \to \cdots$$

for the central extension. The value of the genus-$g$ BSCC map under the $d \oplus \tau_{g,1}(1)$ is $(4g(g-1), 0) \in 8\mathbb{Z} \oplus U$. Therefore, the image

$$\mathrm{Im}(H^1(\mathrm{Ch}_{g,1};\mathbb{Z}) \cong 8\mathbb{Z} \oplus U \to H^1(\mathbb{Z};\mathbb{Z}))$$

of the natural homomorphism is generated by the homomorphism defined by $1 \mapsto \frac{4g(g-1)}{8} = \frac{g(g-1)}{2}$, and the cokernel is isomorphic to the cyclic group of order $\frac{g(g-1)}{2}$, and the Euler class $e \in H^2(\mathrm{Ch}_{g,*};\mathbb{Z})$ is a $\frac{g(g-1)}{2}$-torsion element in the second integral cohomology group $H^2(\mathrm{Ch}_{g,*};\mathbb{Z})$. □

Applying the universal coefficient theorem to the preceding, we obtain the following corollary.

**Corollary 6.14** *For $g \geq 6$, the abelianization $H_1(\mathrm{Ch}_{g,*};\mathbb{Z}) \cong (\mathrm{Ch}_{g,*})^{ab}$ of the Chillingworth subgroup $\mathrm{Ch}_{g,*}$ for the genus-$g$ surface with a base point has $\frac{g(g-1)}{2}$-torsion elements.*

## Acknowledgements

The author would like to thank Takuya Sakasai and Quentin Faes for their valuable discussions and constructive feedback. The author also appreciates Nariya Kawazumi for his useful comments. The author is also grateful to the anonymous referee for a careful review and suggestions, which have helped to improve the clarity and quality of this paper.

RYOTARO KOSUGE kosuge.ryotaro@gmail.com
*Graduate School of Mathematical Sciences, University of Tokyo, Tokyo, Japan*